\PassOptionsToPackage{dvipsnames}{xcolor}

\documentclass[11pt,authoryear,elsevier,numafflabel]{elsarticle}

\usepackage{geometry}
\usepackage{ulem}
\usepackage{xcolor}
\usepackage{setspace}
\usepackage{pdflscape}
\usepackage[round]{natbib}
\usepackage{amsmath}
\usepackage{mathtools}
\usepackage{mathrsfs}
\usepackage{amssymb}
\usepackage{bbold}
\usepackage{colortbl}
\usepackage{booktabs}
\usepackage{multirow}
\usepackage{multicol}

\DeclareMathOperator*{\argmin}{argmin}

\usepackage{amsthm}
\theoremstyle{definition}
\newtheorem{definition}{Definition}

\theoremstyle{remark}
\newtheorem*{remark}{Remark}

\usepackage{algpseudocode}
\usepackage{algorithm}

\usepackage[english]{babel}
\usepackage{csquotes}
\usepackage[T1]{fontenc}
\usepackage[colorlinks=true, allcolors=blue]{hyperref}
\usepackage{subcaption}
\usepackage{fancyhdr}
\usepackage[toc,page]{appendix}
\usepackage{float}

\usepackage{parskip}
\usepackage{a4wide}

\usepackage{tikz}
\usetikzlibrary[arrows.meta]
\usetikzlibrary {decorations.pathmorphing,decorations.pathreplacing,decorations.shapes,decorations.text,shadows,math}
\usetikzlibrary {graphs,positioning,shapes.geometric,quotes}
\usetikzlibrary {datavisualization.formats.functions} 

\title{Heuristic approaches for solving a bilevel optimistic scheduling problem on parallel machines}
\date{}

\author[aff1,aff2]{Quentin Schau\corref{cor1}}
\author[aff2,aff3]{Federico Della Croce}
\author[aff1]{Olivier Ploton}
\author[aff1]{Vincent T'kindt}

\address[aff1]{University of Tours, LIFAT, France.}
\address[aff2]{Politecnico di Torino, DIGEP, Torino, Italy.}
\address[aff3]{CNR, IEIIT, Torino, Italy}
\cortext[cor1]{Corresponding author. \\ E-mail addresses: quentin.schau@univ-tours.fr; quentin.schau@polito.it (Q. Schau)}
\begin{document}

\begin{abstract}
	This work addresses the uniform parallel machine scheduling problem within an optimistic bilevel optimization framework. The leader seeks to minimize the weighted number of tardy jobs, while the follower aims to minimize the total completion time across a set of uniform machines. The hierarchical decision-making process of the bilevel problem makes designing effective heuristics challenging. To tackle this, we exploit a property of the follower that enables the construction of optimal schedules. From this property, we derive an effective branching scheme that simultaneously accounts for both leader and follower decisions. This branching scheme allows us to design a Recovering Beam Search (RBS), which represents a significant contribution from a bilevel perspective. Then we propose a Multi-Start Local Search (MSLS) algorithm based on an innovative scheme that couples the RBS and a Local Search (LS). To the best of our knowledge, while hybridizing beam search with local search is known in other contexts, our approach of leveraging a bilevel-specific branching scheme to efficiently explore the search space is novel. Moreover, we propose an automated method for determining heuristic parameters via Bayesian optimization. This reduces computational resource requirements and yields better parameters than the usual empirical approach. Finally, computational experiments are presented for instances with up to 500 jobs and 10 machines.
\end{abstract}


\def\labelenumi{({\roman{enumi}})}

\maketitle

\section{Introduction}\textit{}
\label{sec:intro}

Bilevel scheduling is a research field at the intersection of scheduling theory and bilevel optimization.
Bilevel scheduling problems involve processing a set of jobs on a set of resources, leading to a set of decisions.
This set of decisions is partitioned between two agents: the leader's set of decisions $X^L$ and the follower's set of decisions $X^F$.
The leader makes the first move by making a decision $x^L \in X^L$, attempting to minimize its own objective function $f^L$.
Next, the follower observes the leader's decision and responds by choosing a decision $x^F \in X^F$ to minimize its own objective function $f^F$.
This results in a hierarchical decision-making process in which both objective functions can only be evaluated after all decisions are taken.

In bilevel optimization problems, the set of optimal solutions to the follower's problem may contain more than one solution.
In this case, the optimal solution computed by the follower may impact the value of the leader's objective function $f^L$.
Thus, to provide a well-defined optimization framework for such problems, we need to make extra assumptions. Several scenarios can be classically considered.
The optimistic case (OPT) corresponds to the scenario where the follower chooses a solution from its set of optimal solutions that is the best for the leader's objective function.
On the contrary, the pessimistic case (PES) corresponds to the scenario where the follower chooses a solution from its set of optimal solutions that is the worst for the leader's objective function. At last, the adversarial case (ADV) corresponds to the scenario where the leader takes decisions to ensure that the optimal solution of the follower's problem is the worst possible, leading to the classical $\min-\max$ problems.

To the best of our knowledge, the literature on algorithms for bilevel scheduling problems is relatively limited. Specifically, exact algorithms are considered in \cite{karlofBilevelProgrammingApplied1996, abassBilevelProgrammingApproach2005a, brownComplexityDelayingAdversarys2005, kovacsConstraintProgrammingApproach2011a, tkindtSingleMachineAdversarial2024,schauBilevelOptimisation2025} and heuristic approaches are considered in \cite{lukacProductionPlanningProblem2008, kisBilevelMachineScheduling2012a, konurAnalysisDifferentApproaches2013a, biancoGridSchedulingBilevel2015}. To solve a bilevel scheduling problem heuristically, we explore the leader's set of decisions. Given a leader's decision, the follower is expected to solve its own problem correspondingly. Note that this can be challenging when the follower's problem is $\cal{NP}$-hard.
If the optimality constraint for the follower's problem is relaxed in the bilevel formulation, then only a lower bound is obtained that, in practice, is often far from the optimal value \citep{colsonOverviewBileveloptimization2007}.

The bilevel scheduling framework has several applications. For instance, unlike usual single-level methods, the bilevel approach enables the modeling of hierarchical decision-making processes within manufacturing systems. In such a system, a leader agent makes initial decisions, such as dispatching production orders, to subordinate follower agents who are then responsible for scheduling those orders. Moreover, we can integrate some maintenance decisions into the scheduling process by requiring that potentially defective machines operate at their minimum speed, while healthy machines work at their maximum speed.

We consider the optimistic bilevel scheduling problem, where the leader selects a subset of jobs and minimizes the weighted number of tardy jobs. Then the follower schedules the selected jobs and minimizes the sum of completion times. The problem can be stated as follows. A set $\cal{J}$ of $N$ jobs is given, and $n\leqslant N$ jobs must be processed on $m_1$ machines with high speed $V_1$ and on $m_0$ machines with low speed $V_0$. Correspondingly, we have $m=m_1+m_0$ machines that can process one job at a time without idle time.
Preemption is not allowed. Each job $J_j$ is characterized by its processing time $p_j$, weight $w_j$, and due date $d_j$. Thus, for a given machine $M_i$ with a speed $V_i$, the processing time of a job $J_j$ on this machine is given by $p_j/V_i$. The follower's objective is to minimize the sum of completion times $\sum_j C^F_j$ where $C^F_j$ denotes the completion time of job $J_j$. The leader's objective is to minimize the weighted number of tardy jobs $\sum_j w_j U^L_j$ where $U^L_j=1$ if job $J_j$ is late, and $U^L_j=0$ otherwise. Let $\mathfrak{S}_{\mathcal{I}}$ be the set of machine schedules that contain jobs from a set of jobs $\mathcal{I}$.
The optimistic bilevel scheduling problem can be formulated as:
\begin{equation} \label{eq:bilevel_problem}
	\begin{aligned}
		 & \min_{\mathcal{I},s^*} \hspace{2mm} \sum_{j \in \mathcal{I}} w_jU^L_j(s^*)                                                                                       \\
		 & \hspace{8mm} \text{s.t.} \quad \mathcal{I} \subset \mathcal{J} \quad \& \quad |\mathcal{I}| = n                                                                  \\
		 & \hspace{8mm} \quad \quad s^* \in \argmin_{s \in \mathfrak{S}_{\mathcal{I}}} \left( Lex \left(\sum_{j \in s} C^F_j(s), \sum_{j \in s} w_j U^L_j(s) \right)\right) \\
	\end{aligned}
\end{equation}

The optimistic scenario implies that there is a lexicographical objective function \\$Lex \left( {\sum_j C^F_j}, {\sum_j w_jU^L_j} \right)$, where we seek a schedule that minimizes the leader's objective function among all optimal schedules for the follower's problem. We refer to this lexicographical problem as the optimistic follower's problem. The bilevel problem is denoted by $Q|OPT-n, V_i \in \{V_{0},V_{1}\}|\sum_j C^F_j,\sum_j w_jU^L_j$, and is known to be $\mathcal{NP}$-hard in the strong sense, and exact algorithms developed in \cite{schauBilevelOptimisation2025} can optimally solve instances with up to 40 jobs for 2 and 4 machines.

	The follower's problem is a special case of the $Q || \sum_j C^F_j$ problem. In \cite{conwayTheoryScheduling1967b}, a characterization of all optimal schedules for the follower's objective function is provided, along with an algorithm to solve it in $\mathcal{O}(n \log n)$ time. These optimal schedules are structured into blocks, where each block consists of a set of pairs (machine index, position on the machine's schedule). Between consecutive blocks, the processing times are increasing. An illustrative example of this block structure is provided in Section \ref{subsec:BlockStruct}.

	The contribution of this work is the presentation of three heuristics designed to tackle this problem. To our knowledge, this represents the first contribution in the bilevel scheduling literature where the optimality of the follower's objective function is guaranteed while the optimistic scenario is relaxed (i.e., the second criterion in the lexicographical problem is not necessarily minimal). This relaxation is motivated by the strong NP-hardness of the optimistic follower's problem and by the structural tractability of the $\sum_j C^F_j$ objective function. Furthermore, these heuristics have a twofold objective: to develop efficient methods for large-scale instances and to discover improved solutions for particularly challenging cases, thereby surpassing the performance of exact algorithms.

	The paper proceeds as follows.
	\hyperref[sec:LS]{Section \ref{sec:LS}} introduces a Local Search (LS) algorithm that attempts to improve a solution by reconsidering the leader and follower decisions, leading to the definition of several neighborhoods.
	\hyperref[sec:RBS]{Section \ref{sec:RBS}} details a Recovering Beam Search (RBS) algorithm, which serves as an improvement over the filtered beam search---a truncated branch-and-bound approach. To the best of our knowledge, this heuristic has not been previously applied to bilevel scheduling problems.
	\hyperref[sec:MSLS]{Section \ref{sec:MSLS}} introduces a Multi-Start Local Search (MSLS) algorithm. This MSLS initially explores the search space using RBS to construct a list of starting solutions and then runs the LS algorithm on each of those solutions.
	\hyperref[sec:ExperRes]{Section \ref{sec:ExperRes}} reports on computational experiments conducted using randomly generated instances. Heuristics' parameterization with Bayesian optimization is also described in this section.
	Finally, \hyperref[sec:CCL]{Section \ref{sec:CCL}} draws conclusions and final remarks.

\section{Local Search}
\label{sec:LS}

\subsection{Block structure example}
\label{subsec:BlockStruct}

In this section, we illustrate the block structure described by \cite{conwayTheoryScheduling1967b} through a concrete example. These blocks are computed by constructing the schedule from right to left, i.e., from the end of the machine schedule towards the beginning. To obtain the final schedule and eliminate idle time, all machine schedules are shifted to the left so that they start at time 0.
We denote such a block as $\mathcal{B}_{b} = \{(i,k)\}$, where $i$ is the machine index and $k$ is the position on the machine schedule. Consider the 9-job example given in Table~\ref{tab:9JobsInstancesExample} to illustrate this block structure. There are two high-speed machines, $M_1$ and $M_2$ with speed $V_1=2$, and two low-speed machines, $M_3$ and $M_4$ with speed $V_0=1$.

\begin{table}[ht!]
	\centering
	\begin{tabular}{cccccccccc}
		\hline
		Job   & $J_1$ & $J_2$ & $J_3$ & $J_4$ & $J_5$ & $J_6$ & $J_7$ & $J_8$ & $J_9$ \\
		\hline
		$p_j$ & 2     & 3     & 4     & 5     & 7     & 8     & 11    & 12    & 13    \\
		$d_j$ & 5     & 1     & 3     & 6     & 4     & 10    & 9     & 13    & 8     \\
		$w_j$ & 1     & 3     & 1     & 1     & 2     & 3     & 2     & 2     & 1     \\
		\hline
	\end{tabular}
	\caption{\color{Green} The example instance consists of 9 jobs}
	\label{tab:9JobsInstancesExample}
\end{table}

Firstly, all jobs are sorted according to the Shortest Processing Time (SPT) rule, i.e., $p_1 \leqslant p_2 \leqslant \hdots \leqslant p_9$. Secondly, starting from the last block $\mathcal{B}_4$, take the $|\mathcal{B}_4|=2$ last jobs ---namely, $J_8$ and $J_9$--- and assign them to this block. Proceed to the previous block $\mathcal{B}_3$ and fill it with the $|\mathcal{B}_3|=4$ last jobs ---$J_7$, $J_6$, $J_5$, and $J_4$. Repeat this process until all blocks are filled. The first block may not be completely filled. The cardinality of all blocks $\mathcal{B}_b$ is computed in $\mathcal{O}(n \log n)$ time in \cite{conwayTheoryScheduling1967b}. Finally, shift the schedule to the left to eliminate any idle time and ensure that it begins at time 0. Figure~\ref{fig:realScheduleExample} shows the final schedule, displaying the completion time $C^F_j$ for each job, and demonstrates that $\sum_j C^F_j = 58$. Figure~\ref{fig:schedulingPatternFProb} illustrates the block structure, where jobs within the dashed group signify their membership to a block and can subsequently undergo permutation without changing the value of $\sum_j C^F_j$.

\begin{figure}[h]
	\hspace{-1cm}
	\begin{minipage}{.5\textwidth}
		\scalebox{0.5}{%
			\begin{tikzpicture}
				\def\yM{-6} 
				\def\sepMachine{0.5} 
				\def\heightMachine{1.25} 
				\foreach \M in {1,...,4}{
						\tikzmath{\yUp = \yM - \sepMachine * (\M -1) - \heightMachine*(\M-1);\yDown = \yM - \sepMachine * (\M -1) - \heightMachine*\M; \yMidle = (\yUp + \yDown)/2; }
						\draw (0,\yUp) -- (0,\yDown) -- (15,\yDown);
						\draw (-0.5,\yMidle) node {$M_{\M}$};
					}

				\tikzmath{\M = 1; \yDown = \yM - \sepMachine * (\M -1) - \heightMachine*\M;}
				\def\C{0} 
				\foreach \name / \P [remember=\C as \Clast (initially 0)] in {1/2,3/4,7/11,8/12}{
						\tikzmath{\p = \P / 2;}
						\draw [draw,fill=YellowOrange!20] (\Clast,\yDown) rectangle (\Clast + \p,\yDown + 1);
						\tikzmath{\C = \Clast + \p; \xJ = \Clast + \p * 0.5; \yJ = \yDown + 0.5;}
						\draw (\xJ,\yJ) node (J\name) {$J_{\name}$};
					}
				\tikzmath{\M = 2; \yDown = \yM - \sepMachine * (\M -1) - \heightMachine*\M;}
				\def\C{0} 
				\foreach \name / \P [remember=\C as \Clast (initially 1)] in {2/3,6/8,9/13}{
						\tikzmath{\p = \P / 2;}
						\draw [draw,fill=YellowOrange!20] (\Clast,\yDown) rectangle (\Clast + \p,\yDown + 1);
						\tikzmath{\C = \Clast + \p; \xJ = \Clast + \p * 0.5; \yJ = \yDown + 0.5;}
						\draw (\xJ,\yJ) node (J\name) {$J_{\name}$};
					}
				\tikzmath{\M = 3; \yDown = \yM - \sepMachine * (\M -1) - \heightMachine*\M;}
				\def\C{0} 
				\foreach \name / \P [remember=\C as \Clast (initially 2.5)] in {5/7}{
						\tikzmath{\p = \P;}
						\draw [draw,fill=YellowOrange!20] (\Clast,\yDown) rectangle (\Clast + \p,\yDown + 1);
						\tikzmath{\C = \Clast + \p; \xJ = \Clast + \p * 0.5; \yJ = \yDown + 0.5;}
						\draw (\xJ,\yJ) node (J\name) {$J_{\name}$};
					}
				\tikzmath{\M = 4; \yDown = \yM - \sepMachine * (\M -1) - \heightMachine*\M;}
				\def\C{0} 
				\foreach \name / \P [remember=\C as \Clast (initially 2.5)] in {4/5}{
						\tikzmath{\p = \P;}
						\draw [draw,fill=YellowOrange!20] (\Clast,\yDown) rectangle (\Clast + \p,\yDown + 1);
						\tikzmath{\C = \Clast + \p; \xJ = \Clast + \p * 0.5; \yJ = \yDown + 0.5;}
						\draw (\xJ,\yJ) node (J\name) {$J_{\name}$};
					}

				\tikzmath{\yDown = \yM - \sepMachine * 3 - \heightMachine*4;}
				\draw [blue,dotted,line width=0.5mm,rounded corners=5pt] (0.1,\yM) -- (0.9,\yM) -- (0.9,\yDown - 0.1) -- (0.1,\yDown - 0.1) -- cycle; 
				\draw [blue,dotted,line width=0.5mm,rounded corners=5pt] (1.1,\yM) -- (2.9,\yM) -- (2.9,\yM - \heightMachine -0.2) -- (2.4,\yM - \heightMachine -0.2) -- (2.4,\yM - \heightMachine*2 - \sepMachine - 0.2) -- (1.1,\yM - \heightMachine *2 - \sepMachine - 0.2) -- cycle; 
				\draw [blue,dotted,line width=0.5mm,rounded corners=5pt] (3.15,\yM) -- (8.4,\yM) -- (8.4,\yM - \heightMachine -0.2) -- (6.4,\yM - \heightMachine -0.2) -- (6.4,\yM - \heightMachine*2 - \sepMachine*2) -- (9.4,\yM - \heightMachine*2 - \sepMachine*2) -- (9.4,\yM - \heightMachine*3 - \sepMachine*3) -- (7.4,\yM - \heightMachine*3 - \sepMachine*3) -- (7.4,\yDown - 0.1) -- (2.65,\yDown - 0.1) -- (2.65,\yM - \heightMachine -\sepMachine) -- (3.15,\yM - \heightMachine -\sepMachine) -- cycle; 
				\draw [blue,dotted,line width=0.5mm,rounded corners=5pt] (8.65,\yM) -- (14.35,\yM) -- (14.35,\yM - \heightMachine - 0.2) -- (12.9,\yM - \heightMachine - 0.2) -- (12.9,\yM - \heightMachine *2 - \sepMachine - 0.2) -- (6.62,\yM - \heightMachine *2 - \sepMachine -0.2) -- (6.62,\yM - \heightMachine -\sepMachine) -- (8.7,\yM - \heightMachine -\sepMachine) -- cycle; 

				\draw[blue] (0.5,\yM + 0.3) node (B1) {$\mathcal{B}_1$};
				\draw[blue] (1.9,\yM + 0.3) node (B2) {$\mathcal{B}_2$};
				\draw[blue] (5.8,\yM + 0.3) node (B3) {$\mathcal{B}_3$};
				\draw[blue] (11.5,\yM + 0.3) node (B4) {$\mathcal{B}_4$};
			\end{tikzpicture}
		}
		\caption{Block structure for the solution of $Q|V_i \in \{V_{0},V_{1}\}|\sum_j C^F_j$}
		\label{fig:schedulingPatternFProb}
	\end{minipage}\hspace{1cm}
	\begin{minipage}{.5\textwidth}
		\scalebox{0.5}{%
			\begin{tikzpicture}
				\def\yM{-6} 
				\def\sepMachine{0.5} 
				\def\heightMachine{1.25} 
				\foreach \M in {1,...,4}{
						\tikzmath{\yUp = \yM - \sepMachine * (\M -1) - \heightMachine*(\M-1);\yDown = \yM - \sepMachine * (\M -1) - \heightMachine*\M; \yMidle = (\yUp + \yDown)/2; }
						\draw (0,\yUp) -- (0,\yDown) -- (15,\yDown);
						\draw (-0.5,\yMidle) node {$M_{\M}$};
					}

				\tikzmath{\M = 1; \yDown = \yM - \sepMachine * (\M -1) - \heightMachine*\M;}
				\def\C{0} 
				\foreach \name / \P [remember=\C as \Clast (initially 0)] in {1/2,3/4,7/11,8/12}{
				\tikzmath{\p = \P / 2;}
				\draw [draw,fill=YellowOrange!20] (\Clast,\yDown) rectangle (\Clast + \p,\yDown + 1);
				\tikzmath{\C = \Clast + \p; \xJ = \Clast + \p * 0.5; \yJ = \yDown + 0.5;}
				\draw (\xJ,\yJ) node (J\name) {$J_{\name}$};
				\draw (\C,\yDown + 1.25 ) node {$C^F_{\name}=\C$};
				}
				\tikzmath{\M = 2; \yDown = \yM - \sepMachine * (\M -1) - \heightMachine*\M;}
				\def\C{0} 
				\foreach \name / \P [remember=\C as \Clast (initially 0)] in {2/3,6/8,9/13}{
				\tikzmath{\p = \P / 2;}
				\draw [draw,fill=YellowOrange!20] (\Clast,\yDown) rectangle (\Clast + \p,\yDown + 1);
				\tikzmath{\C = \Clast + \p; \xJ = \Clast + \p * 0.5; \yJ = \yDown + 0.5;}
				\draw (\xJ,\yJ) node (J\name) {$J_{\name}$};
				\draw (\C,\yDown + 1.25 ) node {$C^F_{\name}=\C$};
				}
				\tikzmath{\M = 3; \yDown = \yM - \sepMachine * (\M -1) - \heightMachine*\M;}
				\def\C{0} 
				\foreach \name / \P [remember=\C as \Clast (initially 0)] in {5/7}{
				\tikzmath{\p = \P;}
				\draw [draw,fill=YellowOrange!20] (\Clast,\yDown) rectangle (\Clast + \p,\yDown + 1);
				\tikzmath{\C = \Clast + \p; \xJ = \Clast + \p * 0.5; \yJ = \yDown + 0.5;}
				\draw (\xJ,\yJ) node (J\name) {$J_{\name}$};
				\draw (\C,\yDown + 1.25 ) node {$C^F_{\name}=\C$};
				}
				\tikzmath{\M = 4; \yDown = \yM - \sepMachine * (\M -1) - \heightMachine*\M;}
				\def\C{0} 
				\foreach \name / \P [remember=\C as \Clast (initially 0)] in {4/5}{
				\tikzmath{\p = \P;}
				\draw [draw,fill=YellowOrange!20] (\Clast,\yDown) rectangle (\Clast + \p,\yDown + 1);
				\tikzmath{\C = \Clast + \p; \xJ = \Clast + \p * 0.5; \yJ = \yDown + 0.5;}
				\draw (\xJ,\yJ) node (J\name) {$J_{\name}$};
				\draw (\C,\yDown + 1.25 ) node {$C^F_{\name}=\C$};
				}

			\end{tikzpicture}
		}
		\caption{Solution of the 9-jobs example of the $Q|V_i \in \{V_{0},V_{1}\}|\sum_j C^F_j$ problem}
		\label{fig:realScheduleExample}
	\end{minipage}%
\end{figure}

\subsection{Heuristic Strategy and Definitions}
\label{subsec:LSAlgo}

To solve an optimistic bilevel minimization problem heuristically---that is, to minimize the leader's objective function $f^L$, subject to the optimistic follower's problem $\text{Lex}(f^F, f^L)$ being minimal--- three theoretical classes of heuristics can be considered.
The first class (1) comprises algorithms that heuristically determine decisions for the leader, subsequently solving the optimistic follower's problem to optimality.
The second class (2) consists of approaches that relax the optimality requirement of the follower's problem $f^F$, while ensuring that the returned solution remains the best possible choice from the leader's perspective.
These approaches are prevalent in the literature. Some recent works study these two classes with certain theoretical guarantees \citep{zareBileveloptimizationInexact2020,shiMixedIntegerBilevel2023,clementeLocalSearchBilevel2025}.
The third class (3) relaxes the requirements of class 1 by imposing optimality only on the first criterion of the lexicographical problem: given a leader's decision, it minimizes only the follower's problem $f^F$. Consequently, this approach yields an upper bound on the leader's objective function.

Solving a bilevel problem is inherently more challenging than solving a single-level problem because both the leader and follower problems must often be solved heuristically. Since the first class of heuristics involves solving the $\text{Lex}\left(\sum_j C^F_j, \sum_j w_j U^L_j\right)$ problem---which is $\mathcal{NP}$-hard \citep{schauBilevelOptimisation2025}---we do not investigate this approach. However, a schedule minimizes the total completion time if and only if it follows the block structure. This property allows us to characterize the set of all optimal schedules for the follower's objective function. Our approach belongs to the third class of heuristics, using this block structure to maintain the optimality of the follower's objective function during the search process. Then we attempt to find a better schedule that optimizes the second criterion, $\sum_j w_j U^L_j$.

Local search algorithms are frequently used to tackle combinatorial optimization problems by iteratively improving a current solution within a predefined neighborhood. In our bilevel context, to improve a current solution, we have two possibilities: change the job selection (the leader's decision) or modify the schedule (the follower's decision). This leads to the definition of two distinct neighborhoods.

\begin{definition}
	Let $(\sigma_s,\mathcal{I}_s,\Omega_s)$ be a solution, where $\sigma_s$ represents the schedule, $\mathcal{I}_s$ is the set of scheduled jobs, and $\Omega_s$ is the set of unselected jobs. This solution is said to be \textit{feasible} if and only if $|\mathcal{I}_s|=n$ and the schedule respects the block structure, i.e., it can be expressed as $\sigma_s = \mathcal{B}_1 \cup \hdots \cup \mathcal{B}_{b_{max}}$. We denote by $C^*(s)=\sum\limits_{j \in \sigma_s}C^F_j(s)$ the optimal value of the follower's objective function and $W(s)=\sum\limits_{j\in \sigma_s}w_j U^L_j(s)$ the value of the leader's objective function.
\end{definition}

There exist job groups such that all jobs within each group have equal processing times. In \cite{schauBilevelOptimisation2025}, an exact polynomial-time algorithm is provided to solve the follower's problem when all jobs have equal processing times, $p_j = p$ for all $j$. We use this algorithm whenever this special case has to be solved.

\subsection{Neighborhoods}

We denote by $S_{|\mathcal{B}_b|}$ the set of all permutations of the jobs in any block $\mathcal{B}_b$. We also denote by $\Pi=\prod\limits_{i=1}^{b_{\max}} S_{|\mathcal{B}_i|}$ the set of all sequences of permutations on the blocks.
Given a feasible solution $s=(\sigma_s,\mathcal{I}_s,\Omega_s)$, to obtain a new schedule $\sigma_{s'}$ from $\sigma_s$, it is enough to apply on each block $\mathcal{B}_b$ of $\sigma_s$, a permutation $\pi_b \in S_{|\mathcal{B}_b|}$. Therefore, for a given set of jobs $\mathcal{I}_s$, we can view $\sigma_s$ as a sequence of permutations $\pi=\left(\pi_1,\hdots,\pi_b,\hdots,\pi_{b_{\max}}\right) \in \Pi$, and $\Pi$ describes the set of all feasible schedules for the follower's problem. Consequently, we define the neighborhood of the follower's decision as a mapping $N_F$ from the set $\Pi$ to the set of all possible subsets $2^{\Pi}$ of $\Pi$, i.e., $N_F: \Pi \to 2^{\Pi}$.

Similarly, let $\Psi=\left\{\mathcal{I} |\, \mathcal{I} \subset \mathcal{J} \,, |\mathcal{I}|=n\right\}$ be the set of all feasible subsets of jobs. We define the neighborhood of the leader's decisions as a mapping $N_L$ from the set $\Psi$ to the set of all possible subsets $2^{\Psi}$ of $\Psi$, i.e., $N_L: \Psi \to 2^{\Psi}$.

Necessarily, for any $\sigma_s \in \Pi$, $N_F(\sigma_s) \subset \Pi$ represents a set of neighbors of $\sigma_s$, and for any $\mathcal{I}_s \in \Psi$, $N_L(\mathcal{I}_s) \subset \Psi$ represents a set of neighbors of $\mathcal{I}_s$.

A locally optimal solution from the follower's perspective is a schedule $\sigma_s$ for which no neighbor in $N_F(\sigma_s)$ yields a better value of the leader's objective function. In other words, for a given schedule $\sigma_s$, if there exists a schedule $\sigma_{s'} \in N_F(\sigma_s)$ such that, for all other schedules $\sigma_{s^{''}} \in N_F(\sigma_s)\setminus \{\sigma_{s^{'}}\}$, $\sum_{j \in \sigma_{s'}} w_j U^L_j(\sigma_{s'}) \leqslant \sum_{j \in \sigma_{s^{''}}} w_j U^L_j(\sigma_{s^{''}})$ holds, then $\sigma_{s'}$ is referred to as an improving decision for the follower, leading to a better solution in the neighborhood of $\sigma_s$.

A locally optimal solution for the leader's problem is a selection $\mathcal{I}_s$ that has no neighbor in $N_L(\mathcal{I}_s)$ with a better value of the leader's objective function. In other words, for a given selection $\mathcal{I}_s$, no other neighboring selection, together with its corresponding schedule, yields a better value of the leader's objective function.
This implies that the optimistic follower's problem must be solved optimally. However, due to the computational hardness of this problem, we instead employ a locally optimal solution to the optimistic follower's problem, leading to the concept of a weakly locally optimal solution for the leader's problem.
The weakness arises from the fact that a selection may exist with an associated schedule that yields a smaller value of the leader's objective function, but such a schedule is not identified due to local optimality.
Let $\tilde{\sigma}_{s}$ denote a locally optimal solution to the optimistic follower's problem corresponding to the leader's decision $\mathcal{I}_s$. Thus, for a given leader's decision $\mathcal{I}_s$, if there exists a selection $\mathcal{I}_{s'} \in N_L(\mathcal{I}_{s})$ such that, for all other selections $\mathcal{I}_{s''} \in N_L(\mathcal{I}_{s})\setminus\{\mathcal{I}_{s'}\}$, $\sum_{j \in \mathcal{I}_{s'} } w_j U^L_j(\tilde{\sigma}_{s'}) \leqslant \sum_{j \in \mathcal{I}_{s''}} w_j U^L_j(\tilde{\sigma}_{s''})$ holds, then $\mathcal{I}_{s'}$ is referred to as a leader's improving decision, leading to a better solution in the neighborhood of $\mathcal{I}_{s}$.



Figures \ref{fig:FollowerLocalOpt} and \ref{fig:LeaderLocalOpt} illustrate the definitions of local optimality and weak local optimality. In both figures, the optimal solution is represented by a star (\begin{tikzpicture} \draw[] (0,0) node[star, star point height=0.05pt,star point ratio=0.3pt,draw, red,fill] {}; \end{tikzpicture}), while the locally optimal solution is indicated by a diamond (\begin{tikzpicture} \draw[] (0,0) node[shape=diamond,inner sep=1.5pt,draw,blue,fill] {}; \end{tikzpicture}) for the optimistic follower's problem and the leader's problem. For the leader's problem, we consider the locally optimal solution of the optimistic follower's problem instead of the global optimum, which restricts the search to the set $N_F(\sigma_s)$. Finally, the weakly locally optimal solution is denoted by a circle (\begin{tikzpicture} \draw[] (0,0) node[shape=circle,inner sep=1.5pt,draw,Green,fill] {}; \end{tikzpicture}).

\begin{figure}[h]
	\begin{minipage}{.5\textwidth}
		\begin{tikzpicture}
			\draw (0,3.2) node {$\sum_j w_j U_j$};
			\draw (5.2,-0.5) node {$\sigma_s$};
			\draw[->] (0,0) -- (5,0);
			\draw[->] (0,0) -- (0,3);
			\fill[fill=blue!20,draw=black] (2.5,1.5) ellipse [x radius=65pt, y radius=40pt];
			\draw (4,2) node {$\Pi$};
			\fill[fill=red!20,draw=black,rotate around={-60:(2.5,1.5)}] (2.5,1.5) ellipse [x radius=35pt, y radius=20pt];
			\draw[] (2.3,2.1) node[node font=\tiny] {$N_F(\sigma_s)$};
			\draw[] (2.5,0.12) node[star, star point height=0.05pt,star point ratio=0.3pt,draw, red,fill] {};
			\draw[] (3,0.4) node[shape=diamond,inner sep=1.5pt,draw,blue,fill] {};
		\end{tikzpicture}
		\caption{Illustration of a locally optimal solution for the optimistic follower's problem}
		\label{fig:FollowerLocalOpt}
	\end{minipage}\hspace{1cm}
	\begin{minipage}{.5\textwidth}
		\begin{tikzpicture}
			\draw (0,3.2) node {$\sum_j w_j U_j$};
			\draw (5.2,-0.5) node {$\mathcal{I}_s$};
			\draw[->] (0,0) -- (5,0);
			\draw[->] (0,0) -- (0,3);
			\fill[fill=blue!20,draw=black] (2.5,1.5) ellipse [x radius=65pt, y radius=40pt];
			\draw (4,2) node {$\Psi$};
			\fill[fill=red!20,draw=black,rotate around={60:(2.5,1.5)}] (2.5,1.5) ellipse [x radius=35pt, y radius=20pt];
			\fill[fill=green!20,draw=black,rotate around={180:(2.5,1.2)}] (2.5,1.2) ellipse [x radius=15pt, y radius=7pt];
			\draw[] (2.7,2.1) node[node font=\tiny] {$N_L(\mathcal{I}_s)$};
			\draw[] (2.5,1.2) node[node font=\tiny] {$N_F(\sigma_s)$};
			\draw[] (2.5,0.12) node[star, star point height=0.05pt,star point ratio=0.3pt,draw, red,fill] {};
			\draw[] (2.1,0.4) node[shape=diamond,inner sep=1.5pt,draw,blue,fill] {};
			\draw[] (2.5,0.95) node[shape=circle,inner sep=1.5pt,draw,Green,fill] {};
		\end{tikzpicture}
		\caption{Illustration of a weak locally optimal solution for the leader's problem}
		\label{fig:LeaderLocalOpt}
	\end{minipage}%
\end{figure}

We can explore our neighborhoods iteratively to find better solutions starting from an initial solution. There are two strategies: \textit{first improvement}, when we stop the exploration at the first improving decision, and \textit{best improvement}, when we search for the best improving decision in a given neighborhood. Preliminary results demonstrate that the best improvement strategy is more efficient for our Local Search algorithms.

\subsection{Exploration}
\label{subsec:ExploLS}
In this section, we describe how we construct $N_F$ and $N_L$ and how these neighborhoods are explored. Given a feasible solution $s$ and its schedule $\sigma_s$, to move to a new schedule $\sigma_{s'}$ in $N_F(\sigma_s)$, we must define a sequence of permutations $\pi_b$ on each block $\mathcal{B}_b$.

\begin{figure}[th!]
	\scalebox{0.8}{%
		\begin{tikzpicture}
			\def\yM{-6} 
			\def\sepMachine{0.5} 
			\def\heightMachine{1.25} 
			\foreach \xM / \jobsMa / \jobMb in {0/{1/1,3/3}/{2/2,4/4},4.5/{2/2,3/3}/{1/1,4/4},9/{1/1,4/4}/{2/2,3/3},13.5/{2/2,4/4}/{1/1,3/3}}{
			\foreach \M in {1,...,2}{
					\tikzmath{\yUp = \yM - \sepMachine * (\M -1) - \heightMachine*(\M-1);\yDown = \yM - \sepMachine * (\M -1) - \heightMachine*\M; \yMidle = (\yUp + \yDown)/2; }
					\draw (\xM,\yUp) -- (\xM,\yDown) -- (\xM+3.2,\yDown);
					\draw (\xM-0.5,\yMidle) node {$M_{\M}$};
				}

			\tikzmath{\M = 1; \yDown = \yM - \sepMachine * (\M -1) - \heightMachine*\M;}
			\def\C{0} 
			\foreach \name / \P [remember=\C as \Clast (initially 0)] in \jobsMa{
				\tikzmath{\p = \P / 2;}
				\draw [draw,fill=YellowOrange!20] (\xM+\Clast,\yDown) rectangle (\xM+\Clast + \p,\yDown + 1);
				\tikzmath{\C = \Clast + \p; \xJ = \Clast + \p * 0.5; \yJ = \yDown + 0.5;}
				\draw (\xM+\xJ,\yJ) node (J\name) {$J_{\name}$};
			}
			\tikzmath{\M = 2; \yDown = \yM - \sepMachine * (\M -1) - \heightMachine*\M;}
			\def\C{0} 
			\foreach \name / \P [remember=\C as \Clast (initially 0)] in \jobMb{
				\tikzmath{\p = \P / 2;}
				\draw [draw,fill=YellowOrange!20] (\xM+\Clast,\yDown) rectangle (\xM+\Clast + \p,\yDown + 1);
				\tikzmath{\C = \Clast + \p; \xJ = \Clast + \p * 0.5; \yJ = \yDown + 0.5;}
				\draw (\xM+\xJ,\yJ) node (J\name) {$J_{\name}$};
			}
			}
			\tikzmath{\yDown = \yM - \sepMachine * 2 - \heightMachine*2;}
			\foreach \xM / \Baa / \Bab / \Bba / \Bbb in {0/1/2/3/4,4.5/2/1/3/4,9/1/2/4/3,13.5/2/1/4/3}{
					\pgfmathsetmacro\result{or(equal(\xM,0),equal(\xM,9)) ? int(1) : int(0)}
					\draw [blue!85,dotted,line width=0.4mm,rounded corners=5pt]
					(\xM + 0.05,\yM) -- (\xM + \Baa/2 - 0.05,\yM)
					\ifnum\result=1
					-- (\xM + \Baa/2 - 0.05,\yM - \heightMachine - \sepMachine) -- (\xM + \Bab/2 - 0.05,\yM - \heightMachine - \sepMachine)
					\else
					-- (\xM + \Baa/2 - 0.05,\yM - \heightMachine - 0.2) -- (\xM + \Bab/2 - 0.05,\yM - \heightMachine - 0.2)
					\fi
					-- (\xM + \Bab/2 - 0.05,\yM - \heightMachine*2 - \sepMachine -0.2)  -- (\xM + 0.05,\yM - \heightMachine*2 - \sepMachine -0.2)
					-- cycle; 
					\draw (\xM + \Bab/4,\yM - \heightMachine *2 - \sepMachine - 0.5) node[blue] {$\mathcal{B}_1$};
					\draw [blue!85,dotted,line width=0.4mm,rounded corners=5pt]
					(\xM + \Baa/2 + 0.05,\yM) -- (\xM + \Baa/2 + \Bba/2 - 0.05,\yM)
					-- (\xM + \Baa/2 + \Bba/2 - 0.05,\yM - \heightMachine - 0.2) -- (\xM + \Bab/2 + \Bbb/2 - 0.05,\yM - \heightMachine - 0.2)
					-- (\xM + \Bab/2 + \Bbb/2  - 0.05,\yM - \heightMachine*2 - \sepMachine -0.2)  -- (\xM + \Bab/2 + 0.05,\yM - \heightMachine*2 - \sepMachine -0.2)

					\ifnum\result=1
						-- (\xM + \Bab/2 + 0.05,\yM - \heightMachine - 0.2) -- (\xM + \Baa/2 + 0.05,\yM - \heightMachine - 0.2)
						\else
						-- (\xM + \Bab/2 + 0.05,\yM - \heightMachine - \sepMachine) -- (\xM + \Baa/2 + 0.05,\yM - \heightMachine - \sepMachine)
					\fi
					-- cycle; 
					\draw (\xM + \Bab/2 + \Bbb/4,\yM - \heightMachine *2 - \sepMachine - 0.5) node[blue] {$\mathcal{B}_2$};
				}

		\end{tikzpicture}
	}
	\caption{Neighborhood of the optimistic follower's problem with 4 jobs and 2 identical-speed machines.}
	\label{fig:neighborFollower}
\end{figure}
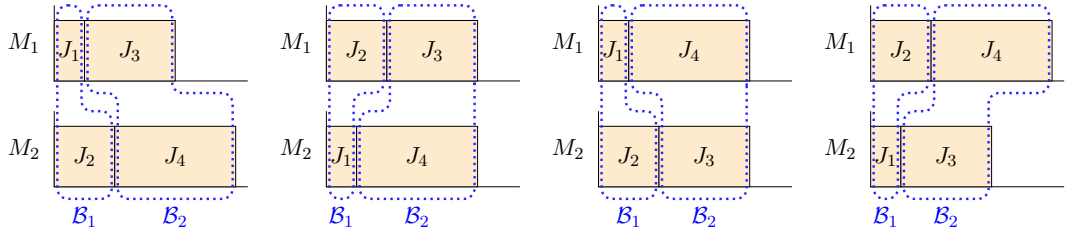

Figure~\ref{fig:neighborFollower} illustrates all the neighbors for a schedule composed of two machines with identical speeds and four assigned jobs. There are only two permutations per block and two blocks, resulting in four distinct machine schedules.

We denote by $\mathcal{J}_b$ the set of jobs scheduled within block $\mathcal{B}_b$ in $\sigma_s$. Let $w_{i,j}$ be the weighted number of tardy jobs on machine $M_i$ when job $J_j$ is scheduled in block $\mathcal{B}_b$, assuming that all other jobs on the machine remain fixed. To compute $w_{i,j}$, we first identify the set $\mathcal{J}_b$ and release the corresponding jobs within block $\mathcal{B}_b$. For each job $J_j \in \mathcal{J}_b$ and each machine $M_i$, we schedule $J_j$ within block $\mathcal{B}_b$ on $M_i$ and update the completion times for the remainder of the schedule on $M_i$. Then we compute the resulting weighted number of tardy jobs on machine $M_i$ as $w_{i,j} = \sum_{J_k \in M_i} w_k U_k(M_i)$. Consequently, our objective is to assign each job $J_j \in \mathcal{J}_b$ to a machine $M_i$ within block $\mathcal{B}_b$ such that the total cost is minimized. This formulation is equivalent to an assignment problem, where the optimal solution yields the permutation $\pi_b$ to be applied to $\mathcal{B}_b$. This first follower's neighborhood is denoted by $N_{F_a}$. The procedure begins with the first block and iteratively solves the assignment problem before proceeding to the next block. Upon reaching the last block, if an improvement is detected, the procedure is re-executed in reverse order.

Because solving the assignment problem has a $\mathcal{O}(n^3)$ time complexity using the Hungarian algorithm \citep{HungarianAlgo}, which can be time expensive, we define a second neighborhood, denoted by $N_{F_s}$, and defined by swapping jobs within each block $\mathcal{B}_b$. In the worst case, we perform $\mathcal{O}(n^2)$ swaps.

From the leader's perspective, we define the neighborhood $N_{L_s}$, by swapping a scheduled job $J_j$ from $\mathcal{I}_s$ with a job $J_r$ from the set of unselected jobs $\Omega_s$. This swap can significantly alter the schedule. Let $\mathcal{B}_{b(j)}$ be the block originally containing job $J_j$, and let $\mathcal{B}_{b(r)}$ be the block where $J_r$ must be scheduled to fulfill the block structure. We address the case where $\mathcal{B}_{b(j)} \neq \mathcal{B}_{b(r)}$ by shifting all jobs. Recall that the block structure imposes that $\max\limits_{J_j \in \mathcal{B}_{b-1}} p_j \leqslant \min\limits_{J_j \in \mathcal{B}_{b}} p_j$ for all $b = 2, \dots, b_{\max}$. Since $\mathcal{B}_{b(r)}$ is full, adding $J_r$ requires shifting one job from $\mathcal{B}_{b(r)}$ to an adjacent block (to the left if $p_j < p_r$ or to the right if $p_j > p_r$). If the adjacent block is also full, continue shifting one job from it to another adjacent block until $\mathcal{B}_{b(j)}$ is reached. When the adjacent block is $\mathcal{B}_{b(j)}$, there is a free position because job $J_j$ has been removed, and no further shifting is necessary.
The job shifted when a block is full is the one with the smallest processing time if it is shifted to the left; otherwise, it is the job with the greatest processing time.

\begin{figure}[th!]
	\centering
	\begin{tikzpicture}
		\def\yM{-6} 
		\def\sepMachine{0.5} 
		\def\heightMachine{1.25} 
		\foreach \xM / \jobsMa / \jobMb in {0/{1/1,6/6}/{2/2,4/4},7/{1/1,6/6}/{4/4,5/5}}{
		\foreach \M in {1,...,2}{
				\tikzmath{\yUp = \yM - \sepMachine * (\M -1) - \heightMachine*(\M-1);\yDown = \yM - \sepMachine * (\M -1) - \heightMachine*\M; \yMidle = (\yUp + \yDown)/2; }
				\draw (\xM,\yUp) -- (\xM,\yDown) -- (\xM+3.2,\yDown);
				\draw (\xM-0.5,\yMidle) node {$M_{\M}$};
			}

		\tikzmath{\M = 1; \yDown = \yM - \sepMachine * (\M -1) - \heightMachine*\M;}
		\def\C{0} 
		\foreach \name / \P [remember=\C as \Clast (initially 0)] in \jobsMa{
			\tikzmath{\p = \P / 2;}
			\draw [draw,fill=YellowOrange!20] (\xM+\Clast,\yDown) rectangle (\xM+\Clast + \p,\yDown + 1);
			\tikzmath{\C = \Clast + \p; \xJ = \Clast + \p * 0.5; \yJ = \yDown + 0.5;}
			\draw (\xM+\xJ,\yJ) node (J\name) {$J_{\name}$};
		}
		\tikzmath{\M = 2; \yDown = \yM - \sepMachine * (\M -1) - \heightMachine*\M;}
		\def\C{0} 
		\foreach \name / \P [remember=\C as \Clast (initially 0)] in \jobMb{
			\tikzmath{\p = \P / 2;}
			\draw [draw,fill=YellowOrange!20] (\xM+\Clast,\yDown) rectangle (\xM+\Clast + \p,\yDown + 1);
			\tikzmath{\C = \Clast + \p; \xJ = \Clast + \p * 0.5; \yJ = \yDown + 0.5;}
			\draw (\xM+\xJ,\yJ) node (J\name) {$J_{\name}$};
		}
		\ifthenelse{\xM=0}{
			\tikzmath{\yDown = \yM - \sepMachine * 2 - \heightMachine*2;}
			\def\sepBlock{0.07}
			\draw [blue!65,dotted,line width=0.4mm,rounded corners=5pt]
			(\xM + \sepBlock,\yM) -- (\xM + 0.5-\sepBlock,\yM)
			-- (\xM + 0.5-\sepBlock,\yM - \heightMachine -\sepMachine) -- (\xM + 1 -\sepBlock,\yM - \heightMachine -\sepMachine)
			-- (\xM + 1 -\sepBlock,\yM - \heightMachine*2 - \sepMachine - 0.2) -- (\xM + \sepBlock,\yM - \heightMachine*2 - \sepMachine - 0.2)
			-- cycle; 
			\draw (\xM + 0.5,\yM - \heightMachine *2 - \sepMachine - 0.5) node[blue] {$\mathcal{B}_1$};
			\draw [blue!65,dotted,line width=0.4mm,rounded corners=5pt]
			(\xM + 0.5 + .05,\yM) -- (\xM + 3.5 -\sepBlock,\yM)
			-- (\xM + 3.5 -\sepBlock,\yM - \heightMachine - 0.2) -- (\xM + 3 - \sepBlock,\yM - \heightMachine - 0.2)
			-- (\xM + 3 -\sepBlock,\yM - \heightMachine * 2 - \sepMachine - 0.2) -- (\xM + 1 + \sepBlock,\yM - \heightMachine * 2 - \sepMachine - 0.2)
			-- (\xM + 1 + \sepBlock,\yM - \heightMachine - 0.2) -- (\xM + 0.5 + \sepBlock,\yM - \heightMachine - 0.2)
			-- cycle; 
			\draw (\xM + 2,\yM - \heightMachine *2 - \sepMachine - 0.5) node[blue] {$\mathcal{B}_2$};
		}{}
		}

		\draw [decoration={brace,mirror,raise=2pt},decorate] (-0.8,\yM) -- (-0.8,\yM - \sepMachine - \heightMachine * 2) node[midway,xshift=-1.5em] {$\sigma_s$};

		\draw (0, \yM + 1) node[] {$\mathcal{I}_s=\left\{J_1,J_2,J_4,J_6\right\}$};
		\draw (4, \yM + 1) node[] {$\Omega_s=\left\{J_3,J_5,J_7\right\}$};
		\draw[line width=1pt,<->,blue] (0.3,\yM + 1.4) .. controls (1,\yM + 2) and (3.5,\yM + 2) .. (4.3,\yM + 1.3) node[pos=0.5, above] {(1)};
		\draw (0.2,\yM+1) node[circle,draw,dashed,inner sep=6pt,blue] {};
		\draw (4.45,\yM+1) node[circle,draw,dashed,inner sep=5pt,blue] {};

		\draw[line width=1pt,->,red] (4.5,\yM + 0.7) .. controls (4.5,\yM-1) and (3.5,\yM - 2) .. (3.2,\yM - \sepMachine - \heightMachine -0.5 ) node[pos=0.5, above,yshift=8pt] {(2)};

		\draw[line width=1pt,<-,red] (0.4,\yM - \sepMachine - \heightMachine -0.2) .. controls (0.5,\yM- \heightMachine) and (2.5,\yM - \heightMachine) .. (2.5,\yM - \sepMachine - \heightMachine -0.1 ) node[pos=0.5, below] {(2)};

		\draw[Green] (5.2, \yM - \sepMachine - \heightMachine) node[label=(3),node font=\Huge] {$\Rightarrow$};
	\end{tikzpicture}
	\caption{Method used to explore the leader's neighborhood.}
	\label{fig:neighborLeader}
\end{figure}
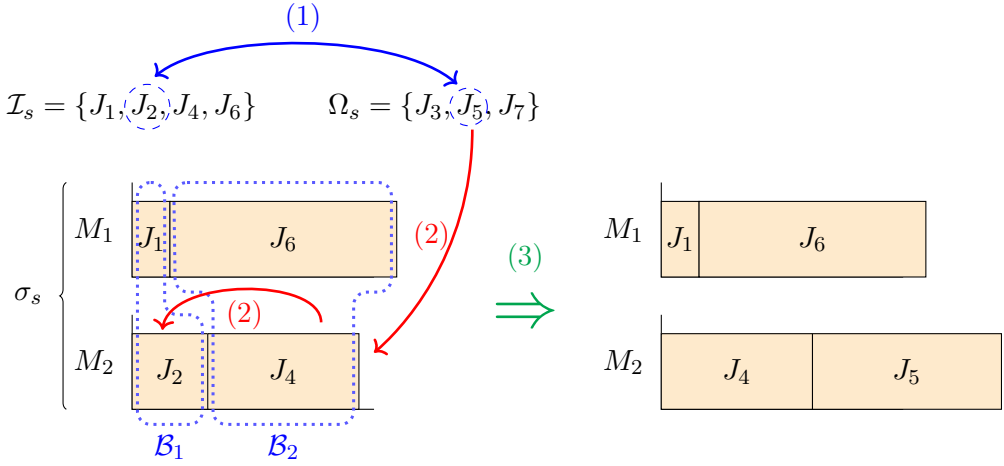

Figure~\ref{fig:neighborLeader} illustrates an example of the method used to explore the neighborhood $N_{L_s}$. We have a list of scheduled jobs $\mathcal{I}_s = \{J_1, J_2, J_4, J_6\}$ and a list of removed jobs $\Omega_s = \{J_3, J_5, J_7\}$. The first step (1) involves swapping a scheduled job, $J_2$, with a removed job, $J_5$. Next, in step (2), we schedule $J_5$ and remove $J_2$. In this example, job $J_5$ must be scheduled in block $\mathcal{B}_2$; however, there is already job $J_4$ present. Thus, we have to shift $J_4$ to $\mathcal{B}_1$ to respect the block structure. Consequently, we obtain the final schedule.

\begin{remark}
	In the literature, neighborhood operators are designed to create neighbors close to the initial solution. However, swapping jobs that are not in the same block can lead to drastic changes in the schedule. Consequently, a more restrictive neighborhood has been considered, where only jobs in the same block are swapped, creating a new schedule with a minimal change---only a swap is performed within the block. However, preliminary experiments have shown that this type of neighborhood does not perform well.
\end{remark}

Our Local Search algorithm (LS)
begins with an initial solution constructed by sorting the jobs in increasing order of the ratio $\frac{p_j-d_j}{w_j}$. The first $n$ jobs from this sorted list are then selected to form an initial set of scheduled jobs. These jobs are scheduled according to the SPT-FAM rule to obtain the first feasible solution. Afterwards, we seek the best weakly locally optimal solution starting from this initial solution by exploring the neighborhood $N_{L_s}$ of the leader's decisions. We consider two versions of the local search, differentiated by the method used to explore the follower's neighborhood: $LS_a$ when the $N_{F_a}$ neighborhood is used, and $LS_s$ when the $N_{F_s}$ neighborhood is used. Due to the potential for an exponential number of improving steps required to converge, the LS algorithm terminates after reaching a predefined maximum number of iterations.

We have also considered a variation of $LS_a$, where the assignment problem is solved not only with jobs already within the block but also by considering all jobs that could be scheduled within the block, particularly those that are currently unselected. This allows for simultaneous decision-making---selecting and scheduling jobs---by solving a single assignment problem. We denote this variation, which solves a "full" assignment problem, as $LS_{fa}$.

\section{Recovering Beam Search}
\label{sec:RBS}

A filtered beam search algorithm, as introduced by \cite{owFilteredBeamSearch1988}, is a truncated branch-and-bound algorithm that explores promising solution sets at each level of a search tree. An improvement to this algorithm, called a Recovering Beam Search (RBS) heuristic, was introduced by \cite{DellaCroceRecoveringBeamSearch2004}. To mitigate the impact of potentially suboptimal selections in the filtered beam search, a recovering mechanism is incorporated. RBS heuristics have been successfully applied to several single-level scheduling problems \citep{dellacroceRecoveringBeamSearch2002,DellaCroceRecoveringBeamSearch2004,valenteFilteredRecoveringBeam2005,ghirardiMakespanMinimizationScheduling2005,esteveRecoveringBeamSearch2006}. However, no such work exists for bilevel scheduling, and the literature remains very limited for bilevel optimization in general. This is explained by the inherent difficulty of managing a branching scheme within a bilevel configuration:
Constraining the leader’s decisions—such as in a branching scheme—does not, in general, yield an a priori bound on the follower’s response, as the follower remains free to optimize subject to its own constraints. Consequently, it is generally difficult to derive strong bounds and, by extension, an effective node selection.

We propose an implementation of the Recovering Beam Search (RBS) tailored to our bilevel problem, using the branching scheme developed by \cite{schauBilevelOptimisation2025}. Each node $\mathcal{N}$ in the search tree represents a partial solution consisting of a partial schedule, a set of selected jobs, a set of removed jobs, and a set of undecided jobs. The branching scheme enforces optimality with respect to the follower's objective function by constructing a schedule that respects the block structure. Consequently, the branching scheme takes leader's and follower's decisions simultaneously: a job is both selected and scheduled at each branching step. The branching procedure involves selecting the next undecided job $J_j$ according to the SPT rule and either assigning $J_j$ to an available position within the first accessible block of the parent node or removing $J_j$, thereby generating a child node $\mathcal{N}_c$. This node $\mathcal{N}_c$ represents all branching decisions. This approach facilitates the construction of a feasible schedule and provides essential state information, such as machine completion times for any given node.

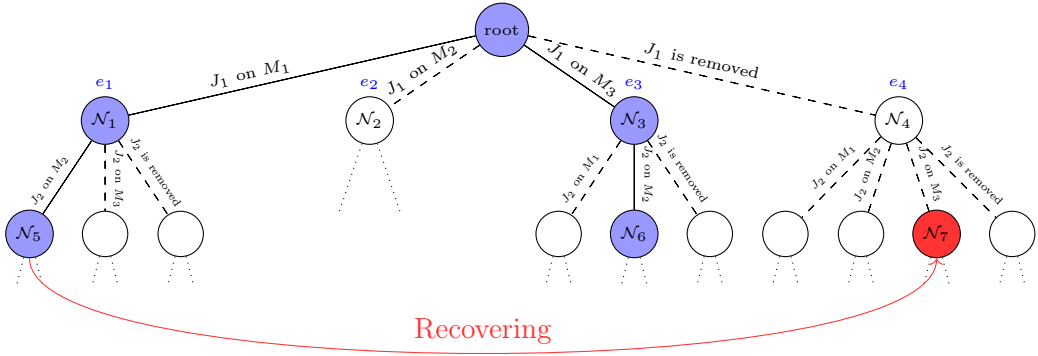
\begin{figure}[ht!]
	\centering
	\begin{tikzpicture}
		\draw
		[every node/.style={circle,draw,solid,node font=\tiny,minimum size=0.6cm},
			level 1/.style={level distance=1.2cm,sibling distance=3.5cm},
			level 2/.style={level distance=1.5cm,sibling distance=1cm},
			level 3/.style={level distance=1cm,sibling distance=0.5cm}
		]
		(0,0) node [fill=blue!40] {root}
		child {
				node [fill=blue!40,label={[blue,yshift=-4pt]$e_1$}] (lvl1-1) {$\mathcal{N}_1$} edge from parent[solid]
				child{
						node [fill=blue!40] (lvl2-1) {$\mathcal{N}_5$} edge from parent[solid]
						child{
								node [draw=none] {} edge from parent[dotted]
							}
						child{
								node [draw=none] {} edge from parent[dotted]
							}
						edge from parent node[right,draw=none,sloped,yshift=5pt,pos=1.05,scale=0.7] {$J_2$ on $M_2$}
					}
				child{
						node [] {} edge from parent[dashed]
						child{
								node [draw=none] {} edge from parent[dotted]
							}
						child{
								node [draw=none] {} edge from parent[dotted]
							}
						edge from parent node[left,draw=none,sloped,yshift=5pt,pos=1.05,scale=0.7] {$J_2$ on $M_3$}
					}
				child{
						node [] {} edge from parent[dashed]
						child{
								node [draw=none] {} edge from parent[dotted]
							}
						child{
								node [draw=none] {} edge from parent[dotted]
							}
						edge from parent node[left,draw=none,sloped,yshift=5pt,pos=1.05,scale=0.7] {$J_2$ is removed}
					}
				edge from parent node[left,draw=none,sloped,yshift=5pt] {$J_1$ on $M_1$}
			}
		child {
				node [label={[blue,yshift=-4pt]$e_2$}] (lvl1-2) {$\mathcal{N}_2$} edge from parent[dashed]
				child {
						node [draw=none] {} edge from parent[dotted]
					}
				child{
						node [draw=none] {} edge from parent[dotted]
					}
				edge from parent node[left,draw=none,sloped,yshift=5pt,pos=0.1] {$J_1$ on $M_2$}
			}
		child {
				node [fill=blue!40,label={[blue,yshift=-4pt]$e_3$}] (lvl1-3) {$\mathcal{N}_3$} edge from parent[solid]
				child{
						node [] {} edge from parent[dashed]
						child{
								node [draw=none] {} edge from parent[dotted]
							}
						child{
								node [draw=none] {} edge from parent[dotted]
							}
						edge from parent node[right,draw=none,sloped,yshift=5pt,pos=1.0,scale=0.7] {$J_2$ on $M_1$}
					}
				child{
						node [fill=blue!40] {$\mathcal{N}_6$} edge from parent[solid]
						child{
								node [draw=none] {} edge from parent[dotted]
							}
						child{
								node [draw=none] {} edge from parent[dotted]
							}
						edge from parent node[left,draw=none,sloped,yshift=5pt,pos=1.,scale=0.7] {$J_2$ on $M_2$}
					}
				child{
						node [] {} edge from parent[dashed]
						child{
								node [draw=none] {} edge from parent[dotted]
							}
						child{
								node [draw=none] {} edge from parent[dotted]
							}
						edge from parent node[left,draw=none,sloped,yshift=5pt,pos=1.,scale=0.7] {$J_2$ is removed}
					}
				edge from parent node[right,draw=none,sloped,yshift=5pt,pos=0.1] {$J_1$ on $M_3$}
			}
		child {
				node [label={[blue,yshift=-4pt]$e_4$}] {$\mathcal{N}_4$} edge from parent[dashed]
				child {
						node [draw] {} edge from parent[dashed]
						child{
								node [draw=none] {} edge from parent[dotted]
							}
						child{
								node [draw=none] {} edge from parent[dotted]
							}
						edge from parent node[right,draw=none,sloped,yshift=5pt,pos=.8,scale=0.7] {$J_2$ on $M_1$}
					}
				child {
						node [draw] {} edge from parent[dashed]
						child{
								node [draw=none] {} edge from parent[dotted]
							}
						child{
								node [draw=none] {} edge from parent[dotted]
							}
						edge from parent node[right,draw=none,sloped,yshift=5pt,pos=1.,scale=0.7] {$J_2$ on $M_2$}
					}
				child {
						node [fill=red!80] (lvl2-2) {$\mathcal{N}_7$} edge from parent[dashed]
						child{
								node [draw=none] {} edge from parent[dotted]
							}
						child{
								node [draw=none] {} edge from parent[dotted]
							}
						edge from parent node[left,draw=none,sloped,yshift=5pt,pos=1.,scale=0.7] {$J_2$ on $M_3$}
					}
				child{
						node [draw] {} edge from parent[dashed]
						child{
								node [draw=none] {} edge from parent[dotted]
							}
						child{
								node [draw=none] {} edge from parent[dotted]
							}
						edge from parent node[left,draw=none,sloped,yshift=5pt,pos=1.,scale=0.7] {$J_2$ is removed}
					}
				edge from parent node[right,draw=none,sloped,yshift=5pt,pos=0.3] {$J_1$ is removed}
			}
		;
		\draw[->,red!80] (lvl2-1) .. controls +(down:2cm) and +(down:2cm) .. node[above,sloped] {Recovering} (lvl2-2);
	\end{tikzpicture}
	\caption{Scheme of the Recovering Beam Search framework}
	\label{fig:RBSProcedure}
\end{figure}

Figure~\ref{fig:RBSProcedure} illustrates the branching scheme for the first two jobs, $J_1$ and $J_2$ in SPT order, for a configuration with three identical machines, $M_1$, $M_2$, and $M_3$. At the first stage, the beam search keeps only node $N_1$, where job $J_1$ is scheduled on $M_1$, and node $N_3$, where $J_1$ is scheduled on $M_3$. Next, at the second stage node $N_5$ (with $J_2$ scheduled on $M_2$) and node $N_6$ (where $J_2$ is also scheduled on $M_2$, although stemming from a different parent) are kept. At any given iteration, child nodes are generated, and only the most promising candidates are preserved for further branching, based on the evaluation phase described in the following section. The recovering mechanism is also detailed in the following section.

\subsection{Evaluation function}

When all child nodes have been created, an evaluation is computed for each node by means of an evaluation function. Then, the RBS retains a subset of $w$ child nodes with the lowest evaluation values at each level of the search tree. This limited number of child nodes, $w$, is referred to as the \textit{beam width} and is a parameter of the algorithm.

For each node $\mathcal{N}$, we denote the associated partial solution by $s_{\mathcal{N}} = (\sigma_{\mathcal{N}}$, $\mathcal{I}_{\mathcal{N}}$, $\Omega_{\mathcal{N}})$, and by $J_{\mathcal{N}}$ the set of undecided jobs. Each generated child node is evaluated by the function $f(\mathcal{N})$ defined as $f(\mathcal{N}) = \alpha LB_{\mathcal{N}} + (1-\alpha) UB_{\mathcal{N}}$, where $\alpha \in [0,1]$ is a given weight that is a parameter of the algorithm. Here, $LB_{\mathcal{N}}$ refers to a lower bound obtained via a column generation algorithm developed by \cite{schauBilevelOptimisation2025}. The columns represent sequences of jobs ordered according to the SPT rule. The master problem consists in selecting the appropriate number of columns to constitute a complete schedule. Consequently, this lower bound relaxes the block structure property by not guaranteeing that processing times increase between successive blocks. It solely enforces an SPT order on each machine and ensures that the required number of jobs is assigned to each machine's schedule. Also, $UB_{\mathcal{N}}$ refers to an upper bound.

We derive a heuristic, $LS_{fac}$, from $LS_{fa}$ that enables the completion of a partial schedule, i.e. a solution with the schedule $\sigma_s$ and $|\mathcal{I}_s| < n$. Using the average processing time, $\bar{p}$, for the set of remaining jobs, we can estimate completion times of all positions in the machine schedule where no job is currently scheduled. Next, we solve an assignment problem within each block, using these estimated completion times, and remove the assigned jobs from the set of remaining jobs. This procedure is repeated from the last to the first unfilled block of the schedule.
Consequently, the upper bound is obtained by completing the partial schedule $\sigma_{\mathcal{N}}$ using the heuristic $LS_{fac}$ to obtain a complementary schedule $\gamma_{\mathcal{N}}$. Thus, a feasible solution is obtained by concatenating the partial schedule $\sigma_{\mathcal{N}}$ and the complementary schedule $\gamma_{\mathcal{N}}$, which we denote by $S_{\mathcal{N}} = \sigma_{\mathcal{N}} || \gamma_{\mathcal{N}}$.

Figure~\ref{fig:RBSProcedure} illustrates the framework of the RBS with a beam width $w=2$. The root node has generated four child nodes: $\mathcal{N}_1$, $\mathcal{N}_2$, $\mathcal{N}_3$, and $\mathcal{N}_4$. These nodes are evaluated using the function $f$, where $e_i = f(\mathcal{N}_i)$ denotes the evaluation of node $\mathcal{N}_i$. Given the computed values $e_3 \leqslant e_1 \leqslant e_2 \leqslant e_4$, nodes $\mathcal{N}_3$ and $\mathcal{N}_1$ are selected as the most promising candidates for the next level of the search tree.

\subsection{Recovering phase}

By selecting only $w$ child nodes, the RBS achieves a suboptimal exploration of the solution space. To recover from potential incorrect selection or evaluation of nodes, a recovering mechanism is introduced. For a given node $\mathcal{N}_1$, among the $w$ selected nodes, the decisions made regarding scheduled jobs $\mathcal{I}_{\mathcal{N}_1}$ and removed jobs $\Omega_{\mathcal{N}_1}$ are reassessed to reach a new node $\mathcal{N}_2$ at the same level of the search tree. A new partial schedule $\sigma_{\mathcal{N}_2}$ is generated by swapping a removed job $J_r \in \Omega_{\mathcal{N}_1}$ with a scheduled job $J_s \in \mathcal{I}_{\mathcal{N}_1}$. This swap necessitates modifying the schedule $\sigma_{\mathcal{N}_1}$ using the same method described for exploring the leader's neighborhood (Figure~\ref{fig:neighborLeader}). The complementary schedule $\gamma_{\mathcal{N}_1}$, previously computed in the evaluation function is then utilized to construct a feasible solution $\sigma_{\mathcal{N}_2} || \gamma_{\mathcal{N}_1}$ and let $UB_{\mathcal{N}_2}$ be the value of this feasible solution. If $UB_{\mathcal{N}_2} < UB_{\mathcal{N}_1}$ then node $\mathcal{N}_2$ is considered for the remainder of the exploration, thus replacing $\mathcal{N}_1$. This recovering process is repeated on the resulting node $\mathcal{N}_2$ until no further improvement is observed.


In Figure~\ref{fig:RBSProcedure}, the recovering phase is illustrated at level 2, enabling a transition from child node $\mathcal{N}_5$ (stemming from $\mathcal{N}_1$) to child node $\mathcal{N}_7$ (stemming from $\mathcal{N}_4$). Consequently, we have reassessed the decision to schedule $J_1$ on $M_1$ by instead removing $J_1$, while also scheduling $J_2$ on $M_3$ rather than $M_2$. This process corresponds to a jump between nodes at the same level of the search tree. Notably, the sibling nodes of $\mathcal{N}_7$ are not generated: specifically, no other child nodes of $\mathcal{N}_4$ are generated. Therefore, for the remainder of the exploration, node $\mathcal{N}_7$ is considered instead of node $\mathcal{N}_5$.

A good set of parameters was identified using Bayesian optimization, as detailed in Section~\ref{subsec:BayeOpt}.


\section{Multi-Start Local Search}
\label{sec:MSLS}

Multi-Start Local Search (MSLS) is an optimization metaheuristic that repeatedly executes a fast local search from various initial solutions, thereby attempting to overcome local optima by exploring diverse regions of the search space. The algorithm operates in two phases: the first phase focuses on an effective diversification by constructing a set of initial solutions to explore new areas of the search space, while the second phase emphasizes intensification to improve these solutions and identify high-quality solutions.

\subsection{Starting solutions}

The first phase consists of populating a list $L$ of at most $K$ initial solutions. The local search algorithm, described in Section~\ref{sec:LS}, finds a good weak locally optimal solution, which is highly dependent on the initial solution. Moreover, the search for a good solution for a local search heuristic is made more complicated by the hierarchical structure in the solution process of the bilevel scheduling problem. Consequently, considering randomly generated initial solutions may not be effective. Instead, we leverage the RBS algorithm, to generate a set of initial solutions. To reduce computational time, we do not use any lower bound during the RBS, but rather prioritize the discovery of diverse solutions.

The list $L$ is populated during the exploration of the search tree by the RBS. We employ two methods to fill $L$: firstly, we include the solutions corresponding to all leaf nodes. Secondly, if the number of leaf nodes $N_{leaf}$ is less than $K$, the list is supplemented with solutions encountered during the evaluation and recovery processes. Specifically, for any node $\mathcal{N}$, multiple solutions may be generated from the evaluation phase and from the recovering mechanism. From these candidate sets, we retain only the $\max(0, K - N_{leaf})$ best-performing solutions, where performance is measured by the number of weighted tardy jobs, $\sum_{J_j \in \sigma_s} w_j U_j(\sigma_s)$. This combination of the branching scheme and the recovering phase ensures an effective exploration of the search space, while retaining only the best solutions guarantees efficient sampling, providing a set of diverse starting points.

These well-separated solutions allow the MSLS to converge toward different local optima during the second phase, ultimately improving overall solution quality.

\subsection{Parameterization of the methods}

Therefore, the MSLS requires three parameters: the number of starting solutions $K$, the beam width $w$, and the ratio of the time budget allocated to the LS. We conducted preliminary experiments with a time budget ratio of 0.5, meaning half of the time is spent on the RBS and the other half in the LS. These preliminary experiments demonstrated that increasing $K$ leads to improved solution quality, although the effect diminishes for larger values of $N$. This suggests that the local search phase contributes positively to the quality of the best solution found by the MSLS. To assess the impact of the beam width $w$, the LS was executed once on the best leaf node solution, i.e., with $K=1$. We observed that increasing $w$ improves solution quality up to a critical value, $w^*$, beyond which the quality deteriorates. Furthermore, the optimal beam width $w^*$ does not allow for a complete exploration of the search space, as the time limit was reached for some instances. Consequently, the beam width $w$ should be maximized as much as possible, up to the point where solution quality deteriorates.

In addition, preliminary experiments show that the computational time and quality are functions of $w$, $K$, and the number of machines, $m$. Moreover, larger values of $w$ and $K$ do not necessarily yield the best results. Consequently, the MSLS can be made more efficient by employing instance-specific parameters. We used Bayesian optimization to determine effective values for these parameters, as described in \hyperref[subsec:BayeOpt]{Section \ref{subsec:BayeOpt}}.

\section{Experimental results}
\label{sec:ExperRes}

In this section, we focus on the experimental evaluation of the heuristic algorithms. All the algorithms were coded in C++, and the tests were run on a PC with a 2 GHz AMD Epyc 7702 processor. The source code can be found on the GitHub repository \citep{schauHeuristicBilevelCode}.
Instances were randomly generated, using the same protocol as described in \citet{framinanOrderSchedulingTardiness2018a} and \citet{hoffmannMinimizingEarlinessTardiness2024}. There are $m_1$ high-speed machines with $V_1 = 2$ and $m_0$ low-speed machines with $V_0 = 1$. The processing times were drawn randomly from the uniform distribution $U[1, 100]$ and weights from the uniform distribution $U[1, 10]$. Due dates were generated using the uniform distribution $U[P \times (1 - tf - rdd/2), P \times (1 - tf + rdd/2)]$ with $tf, rdd \in \{0.2, 0.4, 0.6, 0.8, 1.0\}$ and $P = \frac{\sum_{j=1}^N p_j}{m_1 \times V_1 + m_0 \times V_0}$. Here, $tf$ and $rdd$ define a class of instances, resulting in 25 classes of instances. For each class, we randomly generated 10 instances of size $N = 40, 50, 60, 70, 80,90,100$. For each value of $N$, we set $n = N/4, N/2, 3N/4$, where values were rounded down if fractional. Moreover, we set the total number of machines to $m = 2, 4,10$, ensuring the same number of high-speed and low-speed machines.

\subsection{Bayesian optimization for Parameter Tuning}
\label{subsec:BayeOpt}

Our heuristic algorithms possess parameters that impact their performance and which value is generally made through empirical experiments. The standard approach consists in finding the best parameters across all instances using preliminary experiments. This is the common way to find the best weight $\alpha$ for RBS heuristics \citep{dellacroceRecoveringBeamSearch2002,DellaCroceRecoveringBeamSearch2004,valenteFilteredRecoveringBeam2005,ghirardiMakespanMinimizationScheduling2005,esteveRecoveringBeamSearch2006}. Nevertheless, there exist other approaches, notably Bayesian optimization, that allow us to find good parameters more effectively and automatically.

Assume that we have a heuristic $H$ with $k$ parameters. We denote by $\vec{\beta}_b$ a vector of initial values for these parameters. Thus, $\vec{\beta}_b$ corresponds to the classical approach used to find parameters. Let $H_{\vec{\beta}_b}$ be the corresponding heuristic, called the \textit{baseline} heuristic, parameterized with these found parameters. The objective is to find the values $\vec{\beta}$ such that the heuristic $H_{\vec{\beta}}$ is more efficient than the baseline $H_{\vec{\beta}_b}$. To do so, we generate a training database $\mathfrak{Base}$ with randomly generated instances for several values of $N$, $n$, and $m$. Let $I$ be an instance from the training database $\mathfrak{Base}$. We denote by $UB(\vec{\beta},I)$ the objective value obtained by $H_{\vec{\beta}}$ using parameters $\vec{\beta}$, by $UB(\vec{\beta}_b,I)$ the objective value obtained by the baseline $H_{\vec{\beta}_b}$, and by $W(I)$ the sum of the $n$ greatest job weights of instance $I$. If the condition $\frac{UB(\vec{\beta}_b,I) - UB(\vec{\beta},I)}{W(I)} > 0$ holds for an instance $I$, it means that the heuristic $H_{\vec{\beta}}$ is better than the baseline $H_{\vec{\beta}_b}$. Therefore, to find a good vector of parameters we solve the following maximization problem:

\begin{equation}\label{eq:bayeOptPB}
	\max_{\vec{\beta} \in A \subset \mathbb{R}^k} f(\vec{\beta})= \sum_{I \in \mathfrak{Base} } \left( \frac{UB(\vec{\beta}_b,I) - UB(\vec{\beta},I)}{W(I)} \right)
\end{equation}

where $A$ is the definition set of all parameters. To solve problem (\ref{eq:bayeOptPB}), we use Bayesian optimization that is a framework for finding the extremum of functions for which we don't know the analytical expression. Furthermore, these objective functions can be non-convex, non-differentiable, or lack of a closed-form expression \citep{brochu2010tutorialbayesianoptimizationexpensive}. We utilize the well-known Upper Confidence Bound (UCB) as the acquisition function to select the next point $\vec{\beta}_t$ for exploration at iteration $t$. The parameters $\beta_i$ can be either real-valued variables or, can be integer-valued variables by leveraging the approach of \citep{Garrido_Merch_n_2020}.

For the rest of the section, we describe how Bayesian optimization is used to tune the best parameters for the MSLS and the RBS heuristics. The limitation of a global parameterization approach is that the same parameter values are used for all instances. However, the parameter values should depend on the instance parameters.
For example, the beam width $w$ is linked to the number of child nodes that can be created in the branching scheme, which is approximately $\mathcal{O}(m)$. Moreover, the depth of the search tree decreases as $N$ decreases. Thus, to explore the maximum number of nodes within a given time budget, the beam width must increase as $N$ decreases. Similarly, the parameter $\alpha$ in the RBS evaluation function is used to weight a lower bound and an upper bound, so the quality of the bounds may depend on the number of jobs $N$ in the instance, the number of jobs $n$ to select, and the number of machines $m$. Lastly, we assume that the ratio of the time budget allocated to the LS in the MSLS as well as the number of starting solutions $K$ remain constant, regardless of the instance size. Preliminary results indicate that the beam width decreases and appears to stabilize for high values of $N$ a behavior that can be modeled by an exponential function.

To make the search space independent of the scale of the instance parameters, we first normalize each parameter $x \in \{N,n,m\}$ as $\widetilde{x}=\frac{x-\min(x)}{\max(x)-\min(x)}$, where $\min(x)$ and $\max(x)$ denote the minimum and maximum values of $x$ over the admissible domain. Thus, $\widetilde{x}\in[0,1]$. Since both the beam width $w$ and the RBS weight $\alpha$ must remain within prescribed intervals, we also use the clipping operator defined by $\operatorname{clip}_{[a,b]}(z)=\min(b,\max(a,z))$.

Therefore, we assume that $w$ is a function defined as:

$w(N,m)=\operatorname{clip}_{[1,100]}\left(\left\lfloor e^{\beta_1\widetilde{N}+\beta_2}+\beta_3\widetilde{m}\right\rfloor+\beta_4\right)$.

The term $e^{\beta_1\widetilde{N}+\beta_2}$ preserves the flexibility of the original exponential model while avoiding a direct dependence on the scale of $N$. Depending on the sign and magnitude of $\beta_1$, the beam width may decrease, increase, or remain nearly constant with respect to $N$. The clipping operation guarantees that the beam width effectively used by the heuristic always belongs to $[1,100]$, allowing saturated plateau regions when the non-clipped value is outside this interval.

We also assume that $\alpha$ is defined as a clipped affine function of the normalized instance parameters: $\alpha(N,n,m)=\operatorname{clip}_{[0,1]}\left(\beta_5\widetilde{N}+\beta_6\widetilde{n}+\beta_7\widetilde{m}+\beta_8\right)$. This normalization gives the terms associated with $N$, $n$, and $m$ comparable magnitudes, while the clipping ensures that the weight used in the RBS evaluation function always remains in the interval $[0,1]$.

Therefore, the expressions before clipping in $w(N,m)$ and $\alpha(N,n,m)$ can model either a monotonic increasing or decreasing function, as suggested by preliminary experiments, or a constant function. Moreover, the clipping allows both parameters to reach plateau regions: for $w$, the beam width can be saturated at either $1$ or $100$, and for $\alpha$, the weight can be saturated at either $0$ or $1$. Consequently, this formulation provides a systematic way to formalize the empirical approach and offers greater potential for identifying better parameter configurations.

To determine $A$, the definition set, we look for values of $\beta_i$ such that the corresponding parameters are consistent with the minimal and maximal values of $N$, $n$, and $m$.

For the parameters defining the beam width, the clipping operation guarantees that the value effectively used by the heuristic belongs to $[1,100]$ for every admissible instance. The upper bound $100$ is arbitrary but reasonable in our setting: larger beam widths induce a substantial computational cost, especially when $N$ increases. Therefore, the bounds imposed on $\beta_1$, $\beta_2$, $\beta_3$, and $\beta_4$ are not used to guarantee feasibility directly, but to define a broad and numerically stable search region.

Since $\widetilde{N}\in[0,1]$, we bound the two parameters of the exponential term by $\log(100)\approx4.61$, which leads to $\beta_1,\beta_2\in[-4.61,4.61]$. These bounds allow the exponential component to become either negligible or sufficiently large to reach the upper plateau induced by the clipping. The coefficient $\beta_3$ models the influence of the number of machines. Since this effect is expected to be less significant than the effect of $N$, we bound its contribution by the scale of the machine parameter, namely $\beta_3\in[-\max(m),\max(m)]$. Finally, the additive term is bounded by $\beta_4\in[\![0,100]\!]$, which allows the model to represent a baseline beam width ranging from the minimal to the maximal admissible value.

For the parameters defining $\alpha$, the clipping operation guarantees that the value effectively used by the heuristic belongs to $[0,1]$ for every admissible instance. Since $N$, $n$, and $m$ are normalized before being combined, these coefficients have comparable effects on the value of $\alpha$. This avoids the situation where the coefficient associated with a large-scale parameter such as $N$ dominates the search. In addition, the clipping formulation allows Bayesian optimization to explore both affine behaviors and saturated behaviors, where the weight remains equal to $0$ or $1$ on part of the instance domain. We impose bounds on $\beta_5$, $\beta_6$, $\beta_7$, and $\beta_8$ only to define a broad and numerically stable search region, i.e. $\beta_5, \beta_6, \beta_7, \beta_8 \in [-1,1]$.

The use of a training database helps mitigate overfitting the best parameters and allows for better generalization to other instances. Thereby, we find the best parameters using the training database and then compare all methods on the separate instances described in Section \ref{subsec:computeRes}. The training database is composed of 7875 instances with $N=40,50,60,70,80,90,100$, $m=2,4,10$ and $n=N/4,N/2,3N/4$. On the training database, we obtained the following results:

For the RBS, the baseline parameters are set as $w=1$, $\alpha=0.5$, and we consider a vector $\vec{\beta} = (\beta_1,\beta_2,\beta_3, \ldots, \beta_8) \in A$, where
$$A=[-4.61,4.61] \times [-4.61,4.61] \times [-10,10] \times [\![0,100]\!] \times [-1,1] \times [-1,1] \times [-1,1] \times [-1,1].$$
Here, $[\![a,b]\!]$ represents the set of integers from $a$ to $b$. For this configuration, we obtained a value $f(\vec{\beta})>0$ which indicates an improvement and yields:
\begin{itemize}
	\item $w(N,m) = \operatorname{clip}_{[1,100]}\left(\lfloor e^{-1.15969 \times \widetilde{N} -1.65353} -10.0 \times \widetilde{m} \rfloor + 5\right)$
	\item $\alpha(N,n,m) = \operatorname{clip}_{[0.0,1.0]}\left( -1.0 \times \widetilde{N} -1.0 \times \widetilde{n} -1.0 \times \widetilde{m} + 1.0 \right)$
\end{itemize}

For the MSLS, the baseline parameters are $w=5$, $K=1000,$ and half of the time budget allocated to the LS. We consider a vector $\vec{\beta} = (\beta_1, \ldots, \beta_6) \in A$, where
$$A=[-4.61,4.61] \times [-4.61,4.61] \times [-10,10] \times [\![0,100]\!] \times [\![1,2000]\!] \times [0.01,0.98].$$
In this setup:
\begin{itemize}
	\item $\beta_5$ represents the number of starting solutions,
	\item $\beta_6$ denotes the proportion of time budget allocated to LS.
\end{itemize}

For this configuration, we obtained a value $f(\vec{\beta})>0$ which indicates an improvement and results in:
\begin{itemize}
	\item $w(N,m) = \operatorname{clip}_{[1,100]}\left(\lfloor e^{-1.36890 \times \widetilde{N} + 0.92932} + 1.38553 \times \widetilde{m} \rfloor + 3\right)$
	\item $\beta_5=1459$ and $\beta_6=0.328$.
\end{itemize}

A total of 53 runs for the MSLS and 52 for the RBS algorithms were performed within the Bayesian optimization process, representing a considerable improvement within the computational cost over a grid search approach where parameters are individually increased in small steps.

\begin{figure}[H]
	\centering
	\resizebox{\linewidth}{!}{%
		\input{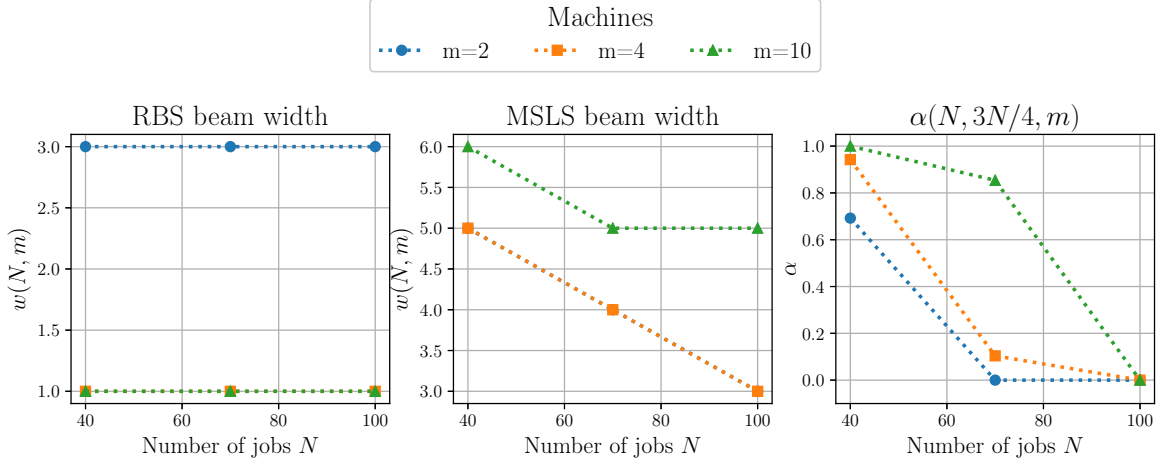}
	}
	\caption{Evolution of the beam width for RBS and MSLS, and of the parameter $\alpha(N,3N/4,m)$, as a function of the number of jobs $N$ for different numbers of machines $m$ and for $n=3N/4$.}
	\label{fig:beam-width-alpha}
\end{figure}

Figure~\ref{fig:beam-width-alpha} shows the evolution of the parameters for $N=40,70,100$, for $m=2,4,10$, and for fixed $n=3N/4$. First, we observe that Bayesian optimization identified, for the RBS, a constant beam width whose value depends on the number of machines. Second, for the MSLS, the beam width decreases for 2 and 4 machines; for 10 machines, it decreases up to 70 jobs, after which it becomes constant and equal to 5. Finally, the parameter $\alpha$ is decreasing: for 2 machines, a plateau equal to 0.0 is reached above 70 jobs, whereas for 4 and 10 machines, the value becomes equal to 0.0 above 100 jobs.

\subsection{Computational results}
\label{subsec:computeRes}

In this section, we present the results of all algorithms. We compare the best exact method, namely the branch-and-bound (BaB) algorithm proposed in \cite{schauBilevelOptimisation2025}, with all heuristic algorithms. All methods were given a time budget of 60 seconds. Since there are many results, we focus only on challenging instances, i.e., instances with $n=3N/4$, and report performance for $N=40,70,100$ and $m=2,4$ and $10$ machines. Tables containing the complete results for additional parameter values are provided in the appendix.

For all tables, we denote by $\Delta = \frac{O^* - O}{O^*}$, where $O^*$ and $O$ represent, respectively, the best-known solution and the current solution being evaluated. Note that, due to the weight of tardy jobs, $\Delta$ can take large values even when the difference between $O^*$ and $O$ is small. For instance, if $O^*=1$ and $O=2$, then $\Delta=50\%$ although the difference is only one unit. To mitigate this effect, we introduce an alternative metric: $\Delta_{\Sigma} = \frac{\sum O^* - \sum O}{\sum O^*}$, where $\sum O^*$ and $\sum O$ represent, respectively, the sum of all best-known objective function values and the sum of the current objective function values. In the tables, the best sum of objective values is presented in bold, and those that outperform truncated BaB are presented in italics.

First, Table~\ref{tab:tabResLS} shows that the LS algorithms are notably fast: in the worst case, the $LS_s$ variant requires 8.41 seconds to solve an instance with $N=100$ jobs and 10 machines, whereas the $LS_{fa}$ variant requires less than 0.1 seconds. Moreover, $LS_s$ slightly outperforms $LS_a$ on instances with 2 machines, with $\Delta_{\Sigma}$ reaching -0.87\%, whereas $LS_a$ tends to outperform $LS_s$ on instances with 4 and 10 machines, with $\Delta_{\Sigma}$ reaching -6.52\%. Thus, in the remainder of this section, we compare $LS_a$ with the other methods.

Note that, when BaB does not reach its time limit, the corresponding instance is optimally solved.

Table~\ref{tab:tabResAll} demonstrates that BaB optimally solves instances with $n=3N/4$ up to $N=40$ jobs.

Second, RBS does not reach its time limit for any value of $N$ on 4 machines. For 2 and 10 machines, the time limit is not reached for values of $N$ up to 70 jobs. Furthermore, for 10 machines, $LS_a$ outperforms RBS in terms of the sum of objective function values. In all other cases, RBS outperforms the LS algorithms.

Finally, combining both heuristics yields the best overall performance. MSLS does not reach its time limit for $N\leqslant70$ with 2 and 4 machines, nor for $N=40$ with 10 machines. When the truncated BaB finds optimal solutions on instances with 2 machines, the gap in the sum of objective function values is less than -0.84\%, the median deviation $\Delta_{med}$ is 0.0\%, and the worst average deviation $\Delta_{avg}$ is -1.22\%. For instances where the truncated BaB outperforms MSLS, i.e., for $N=40$ on 10 machines, while noting that some instances may not be proven optimal because BaB reaches its time limit, MSLS exhibits a gap in the sum of objective function values of less than -0.02\% in the worst case, a $\Delta_{med}$ of 0.0\%, and a $\Delta_{avg}$ of -1.27\%. Consequently, MSLS demonstrates high solution quality both on optimally solved instances and on instances where the truncated BaB finds solutions without a definitive proof of optimality.

For all instances, MSLS outperforms all other heuristics. When MSLS outperforms the truncated BaB, the gap in the sum of objective function values for 2 machines ranges from -26.5\% to -148.55\%. For 4 machines, it ranges from -0.34\% to -122.90\%. For 10 machines, it ranges from -8.55\% to -15.29\%.

The computational times required by the MSLS and RBS algorithms pose a challenge for scaling these methods effectively. However, the LS algorithms, particularly $LS_{fa}$, are designed to handle large instances. Table~\ref{tab:tabResLSHugeInstance} presents computational results for large instances with $N=200, 300, 400, 500$ jobs and $n = 3N/4$, using 2, 4, and 10 machines. We compare the performance of MSLS with that of all local search algorithms. The only cases in which $LS_{fa}$ is outperformed are for $N=200$ on 2 machines, where $LS_s$ performs better, and for $N=200$ on 4 machines, where MSLS performs better. For the other outperformed methods, for $N=200$, the gap in the sum of objective values ranges from -0.09\% to -12.81\% for 2 machines and from -3.36\% to -3.86\% for 4 machines. Moreover, for $N=200$ and 2 machines, $LS_{fa}$ has a better $\Delta_{avg}$ (-25.32\% versus -90.81\% for $LS_s$), while $\Delta_{med}$ is equal to -13.64\% versus 0.00\%. For 4 machines, $LS_{fa}$ has $\Delta_{avg}=-18.25\%$ versus -26.66\%, and $\Delta_{med}=-12.68\%$ versus -8.25\% for MSLS. This means that, when $LS_{fa}$ is not the best method, its average deviation from the best-known solution remains smaller, providing more stable performance. Furthermore, $LS_{fa}$ remains extremely fast, with a worst-case execution time of 0.11 seconds for $N=500$ jobs on 4 machines. Thus, $LS_{fa}$ provides a practical means of solving very large instances with thousands of jobs and tens of machines with minimal computational effort. However, when solution quality is paramount, MSLS should be considered with a time budget scaled to the instance size.

\begin{landscape}
	\begin{table}[ht!]
		\centering
		\fontsize{8}{10}\selectfont
		\setlength{\tabcolsep}{3pt}
		\begin{minipage}{.45\linewidth}
			\begin{tabular}{|c|c|c|cc|c|c|c|c|}
				$N$                    & $n$                   & Method    & $T_{avg}$ & $T_{\max}$ & $\Delta_{avg}$ & $\Delta_{med}$ & $\Delta \sum$ & $\sum$ obj     \\
				\cline{1-9}
				\multicolumn{9}{|c|}{2 Machines}                                                                                                                       \\
				\cline{1-9}
				\multirow[]{3}{*}{40}  & \multirow[]{3}{*}{30} & $LS_a$    & 0.03      & 0.08       & -30.99         & 0.00           & -0.16         & 6968           \\

										&                       & $LS_s$    & 0.02      & 0.04       & -30.74         & 0.00           & 0.00          & \textbf{6957}  \\

										&                       & $LS_{fa}$ & 0.00      & 0.00       & -38.35         & -22.58         & -18.00        & 8209           \\

				\cline{1-9}

				\multirow[]{3}{*}{70}  & \multirow[]{3}{*}{52} & $LS_a$    & 0.36      & 0.83       & -46.87         & -0.82          & -0.87         & 11329          \\

										&                       & $LS_s$    & 0.19      & 0.41       & -46.68         & 0.00           & 0.00          & \textbf{11231} \\

										&                       & $LS_{fa}$ & 0.00      & 0.00       & -33.95         & -21.71         & -15.21        & 12939          \\

				\cline{1-9}

				\multirow[]{3}{*}{100} & \multirow[]{3}{*}{75} & $LS_a$    & 1.65      & 3.53       & -54.77         & 0.00           & -0.18         & 16430          \\

										&                       & $LS_s$    & 0.83      & 1.73       & -51.62         & 0.00           & 0.00          & \textbf{16401} \\

										&                       & $LS_{fa}$ & 0.00      & 0.01       & -34.44         & -21.05         & -11.27        & 18250          \\

				\cline{1-9}
				\multicolumn{9}{|c|}{4 Machines}                                                                                                                       \\
				\cline{1-9}
				\multirow[]{3}{*}{40}  & \multirow[]{3}{*}{30} & $LS_a$    & 0.03      & 0.07       & -18.34         & -1.52          & 0.00          & \textbf{7124}  \\

										&                       & $LS_s$    & 0.02      & 0.05       & -26.78         & -2.94          & -0.93         & 7190           \\

										&                       & $LS_{fa}$ & 0.00      & 0.00       & -34.60         & -20.45         & -16.58        & 8305           \\

				\cline{1-9}

				\multirow[]{3}{*}{70}  & \multirow[]{3}{*}{52} & $LS_a$    & 0.30      & 0.67       & -57.24         & -4.23          & 0.00          & \textbf{11445} \\

										&                       & $LS_s$    & 0.24      & 0.54       & -60.82         & -5.33          & -0.13         & 11460          \\

										&                       & $LS_{fa}$ & 0.00      & 0.00       & -38.08         & -23.75         & -12.36        & 12860          \\

				\cline{1-9}

				\multirow[]{3}{*}{100} & \multirow[]{3}{*}{75} & $LS_a$    & 1.36      & 3.42       & -106.36        & -6.86          & 0.00          & \textbf{16579} \\

										&                       & $LS_s$    & 1.11      & 2.80       & -105.03        & -6.24          & -0.38         & 16642          \\

										&                       & $LS_{fa}$ & 0.00      & 0.00       & -37.96         & -21.97         & -7.46         & 17815          \\

				\cline{1-9}
				\multicolumn{9}{|c|}{10 Machines}                                                                                                                      \\
				\cline{1-9}
				\multirow[]{3}{*}{40}  & \multirow[]{3}{*}{30} & $LS_a$    & 0.06      & 0.11       & -10.76         & 0.00           & 0.00          & \textbf{8614}  \\

										&                       & $LS_s$    & 0.10      & 0.20       & -26.05         & -6.67          & -6.52         & 9176           \\

										&                       & $LS_{fa}$ & 0.00      & 0.00       & -26.62         & -15.69         & -11.60        & 9613           \\

				\cline{1-9}

				\multirow[]{3}{*}{70}  & \multirow[]{3}{*}{52} & $LS_a$    & 0.50      & 1.11       & -23.36         & -4.69          & 0.00          & \textbf{12707} \\

										&                       & $LS_s$    & 0.94      & 2.62       & -43.58         & -9.46          & -6.16         & 13490          \\

										&                       & $LS_{fa}$ & 0.00      & 0.00       & -26.09         & -16.13         & -12.55        & 14302          \\

				\cline{1-9}

				\multirow[]{3}{*}{100} & \multirow[]{3}{*}{75} & $LS_a$    & 1.88      & 4.37       & -65.32         & -4.52          & 0.00          & \textbf{18069} \\

										&                       & $LS_s$    & 3.16      & 8.41       & -108.37        & -12.50         & -6.39         & 19224          \\

										&                       & $LS_{fa}$ & 0.00      & 0.00       & -26.37         & -14.40         & -8.51         & 19607          \\

				\cline{1-9}
			\end{tabular}
			\caption{Comparison of all LS algorithms for $N=40,70,100$ jobs with $n=3N/4$ on 2,4 and 10 machines.}
			\label{tab:tabResLS}
		\end{minipage}
		\hspace{1cm}
		\begin{minipage}{.45\linewidth}
			\begin{tabular}{|c|c|c|cc|c|c|c|c|}
				$N$                    & $n$                   & Method & $T_{avg}$ & $T_{\max}$ & $\Delta_{avg}$ & $\Delta_{med}$ & $\Delta \sum$ & $\sum$ obj     \\
				\cline{1-9}
				\multicolumn{9}{|c|}{2 Machines}                                                                                                                    \\
				\cline{1-9}
				\multirow[]{3}{*}{40}  & \multirow[]{4}{*}{30} & BaB    & 11.64     & 53.00      & 0.00           & 0.00           & 0.00          & \textbf{6580}  \\

										&                       & $LS_a$ & 0.03      & 0.08       & -30.99         & 0.00           & -5.90         & 6968           \\

										&                       & RBS    & 3.58      & 6.95       & -4.66          & 0.00           & -0.94         & 6642           \\

										&                       & MSLS   & 13.18     & 20.89      & -1.22          & 0.00           & -0.84         & 6635           \\

				\cline{1-9}

				\multirow[]{3}{*}{70}  & \multirow[]{4}{*}{52} & BaB    & 49.67     & 60.17      & -62.19         & -5.95          & -26.50        & 13367          \\

										&                       & $LS_a$ & 0.36      & 0.83       & -46.87         & -0.82          & -7.21         & \textit{11329} \\

										&                       & RBS    & 10.50     & 18.66      & -6.26          & 0.00           & -1.14         & \textit{10687} \\

										&                       & MSLS   & 26.61     & 34.80      & -1.72          & 0.00           & 0.00          & \textbf{10567} \\

				\cline{1-9}

				\multirow[]{3}{*}{100} & \multirow[]{4}{*}{75} & BaB    & 52.44     & 60.46      & -330.75        & -174.78        & -148.55       & 38311          \\

										&                       & $LS_a$ & 1.65      & 3.53       & -54.77         & 0.00           & -6.59         & \textit{16430} \\

										&                       & RBS    & 49.08     & 60.01      & -23.35         & -5.26          & -5.64         & \textit{16283} \\

										&                       & MSLS   & 51.77     & 60.02      & -4.97          & 0.00           & 0.00          & \textbf{15414} \\

				\cline{1-9}
				\multicolumn{9}{|c|}{4 Machines}                                                                                                                    \\
				\cline{1-9}
				\multirow[]{3}{*}{40}  & \multirow[]{4}{*}{30} & BaB    & 37.43     & 60.01      & -3.47          & 0.00           & -0.34         & 6816           \\

										&                       & $LS_a$ & 0.03      & 0.07       & -18.34         & -1.52          & -4.87         & 7124           \\

										&                       & RBS    & 0.97      & 2.12       & -10.94         & 0.00           & -2.77         & 6981           \\

										&                       & MSLS   & 11.89     & 20.89      & -3.27          & 0.00           & 0.00          & \textbf{6793}  \\

				\cline{1-9}

				\multirow[]{3}{*}{70}  & \multirow[]{4}{*}{52} & BaB    & 48.02     & 60.07      & -64.12         & -5.88          & -8.73         & 11561          \\

										&                       & $LS_a$ & 0.30      & 0.67       & -57.24         & -4.23          & -7.64         & \textit{11445} \\

										&                       & RBS    & 18.29     & 41.05      & -12.82         & 0.00           & -1.43         & \textit{10785} \\

										&                       & MSLS   & 27.85     & 36.92      & -1.88          & 0.00           & 0.00          & \textbf{10633} \\

				\cline{1-9}

				\multirow[]{3}{*}{100} & \multirow[]{4}{*}{75} & BaB    & 50.80     & 60.33      & -264.16        & -151.91        & -122.90       & 33353          \\

										&                       & $LS_a$ & 1.36      & 3.42       & -106.36        & -6.86          & -10.80        & \textit{16579} \\

										&                       & RBS    & 29.16     & 47.89      & -7.99          & -0.46          & -1.61         & \textit{15204} \\

										&                       & MSLS   & 50.51     & 60.03      & -1.47          & 0.00           & 0.00          & \textbf{14963} \\
				\cline{1-9}
				\multicolumn{9}{|c|}{10 Machines}                                                                                                                   \\
				\cline{1-9}
				\multirow[]{3}{*}{40}  & \multirow[]{4}{*}{30} & BaB    & 32.08     & 60.01      & -2.09          & 0.00           & 0.00          & \textbf{8349}  \\

										&                       & $LS_a$ & 0.06      & 0.11       & -10.76         & 0.00           & -3.17         & 8614           \\

										&                       & RBS    & 0.68      & 1.23       & -16.33         & -4.55          & -5.58         & 8815           \\

										&                       & MSLS   & 18.45     & 24.65      & -1.27          & 0.00           & -0.02         & 8351           \\

				\cline{1-9}

				\multirow[]{3}{*}{70}  & \multirow[]{4}{*}{52} & BaB    & 45.46     & 60.02      & -38.40         & -8.75          & -8.55         & 13231          \\

										&                       & $LS_a$ & 0.50      & 1.11       & -23.36         & -4.69          & -4.25         & \textit{12707} \\

										&                       & RBS    & 27.59     & 47.68      & -26.06         & -5.47          & -4.63         & \textit{12753} \\

										&                       & MSLS   & 52.04     & 60.02      & -3.35          & 0.00           & 0.00          & \textbf{12189} \\

				\cline{1-9}

				\multirow[]{3}{*}{100} & \multirow[]{4}{*}{75} & BaB    & 49.60     & 60.03      & -87.60         & -15.74         & -15.29        & 19675          \\

										&                       & $LS_a$ & 1.88      & 4.37       & -65.32         & -4.52          & -5.88         & \textit{18069} \\

										&                       & RBS    & 50.01     & 60.01      & -50.56         & -13.04         & -9.53         & \textit{18692} \\

										&                       & MSLS   & 53.38     & 60.02      & -11.88         & 0.00           & 0.00          & \textbf{17066} \\
				\cline{1-9}
			\end{tabular}
			\caption{Comparison of BaB algorithm and all heuristic algorithms for $N=40,70,100$ jobs with $n=3N/4$ on 2,4 and 10 machines.}
			\label{tab:tabResAll}
		\end{minipage}
	\end{table}
\end{landscape}

\begin{table}[H]
	\centering
	\fontsize{8}{10}\selectfont
	\begin{tabular}{|c|c|c|cc|c|c|c|c|}
		$N$                    & $n$                    & Method    & $T_{avg}$ & $T_{\max}$ & $\Delta_{avg}$ & $\Delta_{med}$ & $\Delta \sum$ & $\sum$ obj      \\
		\hline
		\cline{1-9}
		\multicolumn{9}{|c|}{2 Machines}                                                                                                                         \\
		\cline{1-9}
		\multirow[]{4}{*}{200} & \multirow[]{4}{*}{150} & $LS_a$    & 39.38     & 60.00      & -86.47         & 0.00           & -0.09         & 32750           \\

								&                        & $LS_s$    & 19.65     & 33.52      & -90.81         & 0.00           & 0.00          & \textbf{32719}  \\

								&                        & $LS_{fa}$ & 0.01      & 0.03       & -25.32         & -13.64         & -7.59         & 35201           \\

								&                        & MSLS      & 57.73     & 60.02      & -62.93         & -17.76         & -12.81        & 36909           \\

		\cline{1-9}

		\multirow[]{4}{*}{300} & \multirow[]{4}{*}{225} & $LS_a$    & 59.88     & 60.01      & -424.54        & -21.83         & -28.79        & 66611           \\

								&                        & $LS_s$    & 58.70     & 60.01      & -302.68        & -9.94          & -14.33        & \textit{59131}  \\

								&                        & $LS_{fa}$ & 0.03      & 0.06       & -1.08          & 0.00           & 0.00          & \textbf{51719}  \\

								&                        & MSLS      & 58.31     & 60.03      & -444.23        & -52.34         & -50.33        & 77748           \\

		\cline{1-9}

		\multirow[]{4}{*}{400} & \multirow[]{4}{*}{300} & $LS_a$    & 60.00     & 60.01      & -518.42        & -40.50         & -43.57        & 99996           \\

								&                        & $LS_s$    & 60.00     & 60.01      & -468.32        & -34.87         & -37.77        & \textit{95956}  \\

								&                        & $LS_{fa}$ & 0.05      & 0.11       & 0.00           & 0.00           & 0.00          & \textbf{69651}  \\

								&                        & MSLS      & 58.90     & 60.03      & -602.84        & -62.47         & -59.98        & 111430          \\

		\cline{1-9}

		\multirow[]{4}{*}{500} & \multirow[]{4}{*}{375} & $LS_a$    & 60.00     & 60.01      & -601.33        & -41.07         & -44.74        & 126595          \\

								&                        & $LS_s$    & 60.00     & 60.01      & -575.78        & -38.10         & -42.15        & \textit{124332} \\

								&                        & $LS_{fa}$ & 0.09      & 0.20       & 0.00           & 0.00           & 0.00          & \textbf{87464}  \\

								&                        & MSLS      & 59.14     & 60.03      & -587.15        & -68.80         & -65.05        & 144363          \\
		\cline{1-9}
		\multicolumn{9}{|c|}{4 Machines}                                                                                                                         \\
		\cline{1-9}
		\multirow[]{4}{*}{200} & \multirow[]{4}{*}{150} & $LS_a$    & 25.41     & 57.86      & -179.15        & -1.66          & -3.92         & 34485           \\

								&                        & $LS_s$    & 21.87     & 49.59      & -144.63        & -1.13          & -3.36         & \textit{34298}  \\

								&                        & $LS_{fa}$ & 0.01      & 0.01       & -18.25         & -12.68         & -3.86         & \textit{34466}  \\

								&                        & MSLS      & 57.09     & 60.02      & -26.66         & -8.25          & 0.00          & \textbf{33184}  \\

		\cline{1-9}

		\multirow[]{4}{*}{300} & \multirow[]{4}{*}{225} & $LS_a$    & 58.61     & 60.01      & -436.82        & -21.19         & -23.97        & 63283           \\

								&                        & $LS_s$    & 59.02     & 60.01      & -442.05        & -19.01         & -22.64        & \textit{62605}  \\

								&                        & $LS_{fa}$ & 0.01      & 0.04       & -0.14          & 0.00           & 0.00          & \textbf{51048}  \\

								&                        & MSLS      & 58.08     & 60.03      & -459.07        & -48.87         & -44.78        & 73909           \\

		\cline{1-9}

		\multirow[]{4}{*}{400} & \multirow[]{4}{*}{300} & $LS_a$    & 59.95     & 60.01      & -582.89        & -38.23         & -41.45        & 96844           \\

								&                        & $LS_s$    & 60.00     & 60.01      & -578.83        & -38.32         & -41.42        & \textit{96824}  \\

								&                        & $LS_{fa}$ & 0.03      & 0.06       & 0.00           & 0.00           & 0.00          & \textbf{68464}  \\

								&                        & MSLS      & 59.15     & 60.03      & -694.85        & -59.57         & -59.81        & 109414          \\

		\cline{1-9}

		\multirow[]{4}{*}{500} & \multirow[]{4}{*}{375} & $LS_a$    & 60.00     & 60.01      & -502.66        & -43.35         & -44.24        & 122638          \\

								&                        & $LS_s$    & 60.00     & 60.01      & -501.24        & -43.35         & -44.18        & \textit{122591} \\

								&                        & $LS_{fa}$ & 0.05      & 0.11       & 0.00           & 0.00           & 0.00          & \textbf{85025}  \\

								&                        & MSLS      & 58.68     & 60.04      & -496.58        & -62.18         & -64.58        & 139933          \\

		\cline{1-9}
		\multicolumn{9}{|c|}{10 Machines}                                                                                                                        \\
		\cline{1-9}
		\multirow[]{4}{*}{200} & \multirow[]{4}{*}{150} & $LS_a$    & 26.90     & 60.01      & -154.63        & -2.43          & -5.74         & 37829           \\

								&                        & $LS_s$    & 38.71     & 60.01      & -182.85        & -11.24         & -12.41        & 40216           \\

								&                        & $LS_{fa}$ & 0.00      & 0.01       & -7.43          & -3.84          & 0.00          & \textbf{35777}  \\

								&                        & MSLS      & 56.91     & 60.02      & -38.33         & -13.61         & -8.52         & 38824           \\

		\cline{1-9}

		\multirow[]{4}{*}{300} & \multirow[]{4}{*}{225} & $LS_a$    & 54.30     & 60.01      & -471.01        & -22.53         & -28.56        & 65795           \\

								&                        & $LS_s$    & 58.03     & 60.01      & -520.47        & -28.37         & -35.11        & 69149           \\

								&                        & $LS_{fa}$ & 0.01      & 0.02       & -0.04          & 0.00           & 0.00          & \textbf{51179}  \\

								&                        & MSLS      & 58.46     & 60.03      & -452.42        & -54.28         & -46.86        & 75159           \\

		\cline{1-9}

		\multirow[]{4}{*}{400} & \multirow[]{4}{*}{300} & $LS_a$    & 58.41     & 60.01      & -459.74        & -35.27         & -38.65        & 93781           \\

								&                        & $LS_s$    & 59.92     & 60.01      & -480.99        & -38.79         & -41.48        & 95697           \\

								&                        & $LS_{fa}$ & 0.01      & 0.03       & 0.00           & 0.00           & 0.00          & \textbf{67638}  \\

								&                        & MSLS      & 59.31     & 60.03      & -509.26        & -65.12         & -59.76        & 108061          \\

		\cline{1-9}

		\multirow[]{4}{*}{500} & \multirow[]{4}{*}{375} & $LS_a$    & 59.83     & 60.01      & -695.31        & -42.30         & -42.72        & 120382          \\

								&                        & $LS_s$    & 60.00     & 60.01      & -719.05        & -43.60         & -44.84        & 122170          \\

								&                        & $LS_{fa}$ & 0.02      & 0.06       & 0.00           & 0.00           & 0.00          & \textbf{84348}  \\

								&                        & MSLS      & 59.08     & 60.04      & -741.53        & -70.68         & -65.81        & 139856          \\

		\cline{1-9}
	\end{tabular}
	\caption{Comparison of all heuristic algorithms for $N=200,300,400$ and 500 jobs with $n=3N/4$ on 2,4 and 10 machines.}
	\label{tab:tabResLSHugeInstance}
\end{table}

\section{Conclusion}
\label{sec:CCL}

We present novel heuristic approaches for solving the optimistic bilevel scheduling problem. Providing polynomial-time heuristics for such problems can be very challenging due to the hierarchical solution process. The heuristics we propose maintain the optimality of the follower's objective function while relaxing the optimistic configuration for the leader. This new perspective on bilevel scheduling theory leverages inherent problem properties, enabling the development of effective heuristics.

A Multi-Start Local Search (MSLS) algorithm uses a Recovering Beam Search (RBS) to efficiently explore the search space with a dedicated branching scheme, and employs a Local Search (LS) algorithm to identify local optima. The MSLS provides high-quality results on instances solved to optimality by an exact branch-and-bound (BaB) algorithm. Furthermore, within a limited time budget, the MSLS outperforms both the BaB and all other tested heuristics on other instances. However, the MSLS becomes slower as the instance size grows, specifically for $N \geq 200$ jobs. In contrast, a Local Search algorithm, $LS_{fa}$, is both fast and scalable, efficiently handling instances with several thousand jobs across tens of machines. Therefore, the MSLS is a quality-oriented method, while $LS_{fa}$ is a scalability-oriented method.

Moreover, we make use of Bayesian optimization to determine good parameters for our heuristics. This approach can be generalized to other heuristics and problems, providing an automated method for fine-tuning parameters. Such an approach can outperform classical methods based on empirical experimentation while requiring a minimal number of evaluations and computational resources. Additionally, it would be interesting to study other $\mathcal{NP}$-hard bilevel scheduling problems and to develop heuristics for tackling them.

\pagebreak

\appendix

\section{Complete Computational Results}

In this section, we provide the complete set of tables reporting the computational results obtained for all combinations of $N=40,50,60,70,80,90,100$, $n=N/4,N/2,3N/4$, and $m=2,4,10$ machines.

\begin{table}[H]
	\centering
	\fontsize{6}{7}\selectfont
	\setlength{\tabcolsep}{0.5pt}
	\makebox[\textwidth][c]{
		\begin{minipage}{.32\linewidth}
			\centering
\begin{tabular}{|c|c|c|cc|c|c|c|c|}
				$N$                    & $n$                   & Method  & $T_{avg}$ & $T_{\max}$ & $\Delta_{avg}$ & $\Delta_{med}$ & $\Delta \sum$ & $\sum$ obj     \\
				\hline

				\multirow[]{9}{*}{40}  & \multirow[]{3}{*}{10} & $LS_a$  & 0.00      & 0.02       & -1.96          & 0.00           & 0.00          & \textbf{207}   \\

									&                       & $LS_s$  & 0.00      & 0.01       & -1.96          & 0.00           & 0.00          & \textbf{207}   \\

									&                       & $LS_{fa}$ & 0.00      & 0.00       & -108.91        & -66.67         & -119.32       & 454            \\

				\cline{2-9}

									& \multirow[]{3}{*}{20} & $LS_a$  & 0.03      & 0.08       & -27.60         & 0.00           & -0.53         & 1713           \\

									&                       & $LS_s$  & 0.01      & 0.03       & -26.39         & 0.00           & 0.00          & \textbf{1704}  \\

									&                       & $LS_{fa}$ & 0.00      & 0.00       & -41.95         & -26.29         & -33.80        & 2280           \\

				\cline{2-9}

									& \multirow[]{3}{*}{30} & $LS_a$  & 0.03      & 0.08       & -30.99         & 0.00           & -0.16         & 6968           \\
\hline
			\end{tabular}
			\caption{Comparison of all LS algorithms for $N=40,50,60,70,80,90,100$ jobs with $n=N/4,N/2,3N/4$ on 2 machines.}
		\end{minipage}
		\hspace{1.6cm}
		\begin{minipage}{.32\linewidth}
			\centering
			\begin{tabular}{|c|c|c|cc|c|c|c|c|}
				$N$                    & $n$                   & Method  & $T_{avg}$ & $T_{\max}$ & $\Delta_{avg}$ & $\Delta_{med}$ & $\Delta \sum$ & $\sum$ obj     \\
				\hline

				\multirow[]{9}{*}{40}  & \multirow[]{3}{*}{10} & $LS_a$  & 0.00      & 0.02       & -10.00         & 0.00           & 0.00          & \textbf{327}   \\

									&                       & $LS_s$  & 0.00      & 0.01       & -11.45         & 0.00           & -0.31         & 328            \\

									&                       & $LS_{fa}$ & 0.00      & 0.00       & -87.10         & -33.33         & -61.77        & 529            \\

				\cline{2-9}

									& \multirow[]{3}{*}{20} & $LS_a$  & 0.02      & 0.07       & -24.94         & 0.00           & 0.00          & \textbf{1947}  \\

									&                       & $LS_s$  & 0.01      & 0.05       & -32.70         & 0.00           & -2.05         & 1987           \\

									&                       & $LS_{fa}$ & 0.00      & 0.00       & -56.73         & -37.20         & -39.39        & 2714           \\

				\cline{2-9}

									& \multirow[]{3}{*}{30} & $LS_a$  & 0.03      & 0.07       & -18.34         & -1.52          & 0.00          & \textbf{7124}  \\
\hline
			\end{tabular}
			\caption{Comparison of all LS algorithms for $N=40,50,60,70,80,90,100$ jobs with $n=N/4,N/2,3N/4$ on 4 machines.}

		\end{minipage}
		\hspace{1.6cm}
		\begin{minipage}{.32\linewidth}
			\centering
			\begin{tabular}{|c|c|c|cc|c|c|c|c|}
				$N$                    & $n$                   & Method  & $T_{avg}$ & $T_{\max}$ & $\Delta_{avg}$ & $\Delta_{med}$ & $\Delta \sum$ & $\sum$ obj     \\
				\hline

				\multirow[]{9}{*}{40}  & \multirow[]{3}{*}{10} & $LS_a$  & 0.01      & 0.03       & -1.14          & 0.00           & 0.00          & \textbf{742}   \\

									&                       & $LS_s$  & 0.02      & 0.06       & -4.28          & 0.00           & -0.94         & 749            \\

									&                       & $LS_{fa}$ & 0.00      & 0.00       & -70.07         & -20.00         & -46.09        & 1084           \\

				\cline{2-9}

									& \multirow[]{3}{*}{20} & $LS_a$  & 0.04      & 0.12       & -7.23          & 0.00           & 0.00          & \textbf{3145}  \\

									&                       & $LS_s$  & 0.06      & 0.20       & -48.71         & -14.81         & -18.63        & 3731           \\

									&                       & $LS_{fa}$ & 0.00      & 0.00       & -39.39         & -17.39         & -23.02        & 3869           \\

				\cline{2-9}

									& \multirow[]{3}{*}{30} & $LS_a$  & 0.06      & 0.11       & -10.76         & 0.00           & 0.00          & \textbf{8614}  \\
\hline
			\end{tabular}
			\caption{Comparison of all LS algorithms for $N=40,50,60,70,80,90,100$ jobs with $n=N/4,N/2,3N/4$ on 10 machines.}

		\end{minipage}
	}
\end{table}

\clearpage

\newgeometry{top=1.5cm,bottom=1.5cm}

\begin{table}[H]
	\centering
	\fontsize{6}{7}\selectfont
	\setlength{\tabcolsep}{0.5pt}
	\makebox[\textwidth][c]{
		\begin{minipage}{.32\linewidth}
			\centering
			\begin{tabular}{|c|c|c|cc|c|c|c|c|}
				$N$                     & $n$                   & Method & $T_{avg}$ & $T_{\max}$ & $\Delta_{avg}$ & $\Delta_{med}$ & $\Delta \sum$ & $\sum$ obj     \\
				\hline

				\multirow[]{12}{*}{40}  & \multirow[]{4}{*}{10} & BaB    & 0.09      & 0.32       & 0.00           & 0.00           & 0.00          & \textbf{203}   \\

										&                       & $LS_a$ & 0.00      & 0.02       & -1.96          & 0.00           & -1.97         & 207            \\

										&                       & RBS    & 0.05      & 0.25       & -0.49          & 0.00           & -0.49         & 204            \\

										&                       & MSLS   & 1.16      & 7.09       & 0.00           & 0.00           & 0.00          & \textbf{203}   \\

				\cline{2-9}

										& \multirow[]{4}{*}{20} & BaB    & 1.62      & 23.20      & 0.00           & 0.00           & 0.00          & \textbf{1543}  \\

										&                       & $LS_a$ & 0.03      & 0.08       & -27.60         & 0.00           & -11.02        & 1713           \\

										&                       & RBS    & 0.92      & 2.36       & -5.56          & 0.00           & -1.62         & 1568           \\
\hline
			\end{tabular}
			\caption{Comparison of BaB algorithm and all heuristic algorithms for $N=40,50,60,70,80,90,100$ jobs with $n=N/4,N/2,3N/4$ on 2 machines.}
		\end{minipage}
		\hspace{1.6cm}
		\begin{minipage}{.32\linewidth}
			\centering
			\begin{tabular}{|c|c|c|cc|c|c|c|c|}
				$N$                     & $n$                   & Method & $T_{avg}$ & $T_{\max}$ & $\Delta_{avg}$ & $\Delta_{med}$ & $\Delta \sum$ & $\sum$ obj     \\
				\hline

				\multirow[]{12}{*}{40}  & \multirow[]{4}{*}{10} & BaB    & 0.10      & 0.21       & 0.00           & 0.00           & 0.00          & \textbf{308}   \\

										&                       & $LS_a$ & 0.00      & 0.02       & -10.00         & 0.00           & -6.17         & 327            \\

										&                       & RBS    & 0.04      & 0.16       & -3.76          & 0.00           & -2.92         & 317            \\

										&                       & MSLS   & 1.55      & 8.98       & 0.00           & 0.00           & 0.00          & \textbf{308}   \\

				\cline{2-9}

										& \multirow[]{4}{*}{20} & BaB    & 2.34      & 45.76      & 0.00           & 0.00           & 0.00          & \textbf{1812}  \\

										&                       & $LS_a$ & 0.02      & 0.07       & -24.94         & 0.00           & -7.45         & 1947           \\

										&                       & RBS    & 0.23      & 0.52       & -9.02          & 0.00           & -3.86         & 1882           \\
\hline
			\end{tabular}
			\caption{Comparison of BaB algorithm and all heuristic algorithms for $N=40,50,60,70,80,90,100$ jobs with $n=N/4,N/2,3N/4$ on 4 machines.}

		\end{minipage}
		\hspace{1.6cm}
		\begin{minipage}{.32\linewidth}
			\centering
			\begin{tabular}{|c|c|c|cc|c|c|c|c|}
				$N$                     & $n$                   & Method & $T_{avg}$ & $T_{\max}$ & $\Delta_{avg}$ & $\Delta_{med}$ & $\Delta \sum$ & $\sum$ obj     \\
				\hline

				\multirow[]{12}{*}{40}  & \multirow[]{4}{*}{10} & BaB    & 0.11      & 0.25       & 0.00           & 0.00           & 0.00          & \textbf{736}   \\

										&                       & $LS_a$ & 0.01      & 0.03       & -1.14          & 0.00           & -0.82         & 742            \\

										&                       & RBS    & 0.05      & 0.19       & -4.89          & 0.00           & -1.36         & 746            \\

										&                       & MSLS   & 2.88      & 16.93      & 0.00           & 0.00           & 0.00          & \textbf{736}   \\

				\cline{2-9}

										& \multirow[]{4}{*}{20} & BaB    & 1.15      & 11.05      & 0.00           & 0.00           & 0.00          & \textbf{3080}  \\

										&                       & $LS_a$ & 0.04      & 0.12       & -7.23          & 0.00           & -2.11         & 3145           \\

										&                       & RBS    & 0.26      & 0.65       & -12.25         & 0.00           & -4.64         & 3223           \\
\hline
			\end{tabular}
			\caption{Comparison of BaB algorithm and all heuristic algorithms for $N=40,50,60,70,80,90,100$ jobs with $n=N/4,N/2,3N/4$ on 10 machines.}

		\end{minipage}
	}
\end{table}

\clearpage
\restoregeometry
\section*{Acknowledgements}

The authors benefited from the use of the cluster at the Centre de Calcul Scientifique en région Centre-Val de Loire.

\bibliographystyle{elsarticle-harv}
\bibliography{./biblio.bib}
\end{document}